\documentclass[12pt]{article}

\usepackage{amsmath,a4wide,hyperref}
\usepackage{amssymb}
\usepackage{amsthm}
\usepackage[english]{babel}
\usepackage[utf8]{inputenc} 

\DeclareUnicodeCharacter{00A0}{~}

\numberwithin{equation}{section}

\newtheorem{Th}{Theorem}
\newtheorem{Lem}{Lemma}
\newtheorem{Prop}{Proposition}

\newtheorem{Remark}{Remark}

\newcommand{\Dem}{\noindent{\bf Proof }}

\makeatletter
\newcommand*{\house}[1]{%
 \mathord{%
 \mathpalette\@house{#1}%
 }%
}
\newcommand*{\@house}[2]{%
 \dimen@=\fontdimen8 %
 \ifx#1\scriptscriptstyle\scriptscriptfont
 \else\ifx#1\scriptstyle\scriptfont
 \else\textfont\fi\fi
 3 %
 \sbox0{%
 $#1%
 \vrule width\dimen@\relax
 \overline{%
 \kern2\dimen@
 \begingroup 
 #2%
 \endgroup
 \kern2\dimen@
 }%
 \vrule width\dimen@\relax
 \mathsurround=1.5\dimen@ 
 $%
 }%
 \ht0=\dimexpr\ht0-\dimen@\relax
 \dp0=\dimexpr\dp0+2\dimen@\relax
 \vbox{%
 \kern\dimen@ 
 \copy0 %
 }%
}

\renewcommand{\P}{\mathbb{P}}
\newcommand{\N}{\mathbb{N}}
\newcommand{\Z}{\mathbb{Z}}
\newcommand{\Q}{\mathbb{Q}}
\newcommand{\Qbar}{\overline{\mathbb{Q}}}
\newcommand{\R}{\mathbb{R}}
\newcommand{\C}{\mathbb{C}}
\newcommand{\K}{\mathbb{K}}
\newcommand{\etoile}{^\ast}

\newcommand{\Span}{{\rm Span}}

\newcommand{\eps}{\varepsilon}

\newcommand{\om}{\omega}
\newcommand{\eneq}{\end{equation}}
\newcommand{\lam}{\lambda}

\newcommand{\unb}{\{1,\ldots,b\}}
\newcommand{\una}{\{1,\ldots,a\}}

\newcommand{\zeron}{\{0,\ldots,n\}}

\newcommand{\calN}{\mathcal{N}}

\newcommand{\ord}{{\rm ord}}
\newcommand{\dd}{{\rm d}}
\newcommand{\rk}{{\rm rk}}

\newcommand{\Card}{{\rm Card}\,}
\newcommand\tra{ \ ^t }

\newcommand{\combi}[2]{{ \left( \begin{array}{c} #1 \\ #2 \end{array} \right)}}

\newcommand{\calC}{{\mathcal{C}}}

\newcommand{\calS}{{\mathcal{S}}}

\newcommand{\vol}{{\rm vol}\, }

\newcommand{\hsoul}{\underline{h}}

\newcommand{\cij}{c_{i,j}}

\newcommand{\del}{\delta}

\newcommand{\lcm}{{\rm lcm}}
\newcommand{\Id}{{\rm Id}}
\newcommand{\frakA}{{\mathfrak A}}

\newcommand{\qq}{q}
\newcommand{\qqp}{Q^{[p]}}
\newcommand{\qqppr}{Q^{[p]_{'}} }
\newcommand{\capa}{\kappa}
\newcommand{\likn}{\ell_{p,k,i}^{(n)}}

\newcommand{\lzkn}{\ell_{p,k,0}^{(n)}}

\newcommand{\lebeta}{\alpha}

\newcommand{\OK}{{\mathcal O}_{\K}}
\newcommand{\OQzz}{{\mathcal O}_{\Q(z_0)}}

\newcommand{\Snpzero}{S_{n,p}^{[0]}}
\newcommand{\Vpzero}{V_p^{[0]}}
\newcommand{\Snpinf}{S_{n,p}^{[\infty]}}
\newcommand{\Vpinf}{V_p^{[\infty]}}
\newcommand{\Vpinfkmu}{V_p^{[\infty] (k-1)}}
\newcommand{\Snpinfkmu}{S_{n,p}^{[\infty] (k-1) }}
\newcommand{\Snpzerokmu}{S_{n,p}^{[0](k-1) }}

\newcommand{\Li}{{\rm Li}}

\newcommand{\zeroh}{\{0,\ldots,h\}}
\newcommand{\intk}{\{2rn+2,\ldots,\capa n\}}
\newcommand{\ff}{f}
\newcommand{\hh}{\varrho}
\newcommand{\yy}{y}
\newcommand{\hhti}{\widetilde \hh}
\newcommand{\yti}{\widetilde y}
\newcommand{\gdyti}{\widetilde Y}
\renewcommand{\SS}{S}
\newcommand{\TT}{T}
\newcommand{\bino}{{q \choose p }} 
\newcommand{\alp}{\alpha}
\newcommand{\qshid}{N}
\newcommand{\cstun}{c_1}

\newcommand{\unlebeta}{\{1,\ldots,\lebeta\}}
\newcommand{\Spe}{\calS_{p,e}}

\newcommand{\bplush}{b+h}
\newcommand{\aplush}{a+h}
\newcommand{\gdom}{\Omega}
\newcommand{\theco}{\vartheta_{k,i,j,i',j'}}
\newcommand{\thecoTIT}{\vartheta_{k,i,T,I,T}}
\newcommand{\thecojIT}{\vartheta_{k,i,j,I,T}}
\newcommand{\thecoze}{\vartheta_{k,0,j,i',j'}}
\newcommand{\thecozeIT}{\vartheta_{k,0,j,I,T}}
\newcommand{\fctchi}{\kappa}
\newcommand{\ateeps}{\alpha_\eps}
\newcommand{\ate}{\alpha_1}
\newcommand{\bte}{\alpha_0}

\title{Linear independence of odd zeta values using Siegel's lemma}
\author{St\'ephane Fischler\footnote{Universit\'e Paris-Saclay, CNRS, Laboratoire de math\'ematiques d'Orsay, 91405 Orsay, France}}
\date{\today}

\begin{document}

\maketitle

\begin{abstract}
We prove that among 1 and the odd zeta values $\zeta(3)$, $\zeta(5)$, \ldots, $\zeta(s)$, at least $ 0.21 \sqrt{s / \log s}$ are linearly independent over the rationals, for any sufficiently large odd integer $s$. This is the first asymptotic improvement on the lower bound, logarithmic in $s$, obtained by Ball-Rivoal in 2001.

The proof is based on Siegel's lemma to construct non-explicit linear forms in values at odd integers of the Riemann zeta function, instead of using explicit well-poised hypergeometric series. A new refinement of Siegel's linear independence criterion  is applied, together with a multiplicity estimate (namely a generalization of Shidlovsky's lemma).

The result is also adapted to deal with values of the first $s$ polylogarithms at a fixed algebraic point in the unit disk, improving bounds of Rivoal and Marcovecchio.
\end{abstract}

\bigskip

\noindent Math. Subject Classification: 11J72 (Primary), 11M06 (Secondary).

\bigskip

\section{Introduction}

It is well known that $\zeta(s) = \sum_{n=1}^\infty n^{-s} $ is equal, when 
 $s\geq 2$ is an even integer, to $c_s \pi^{ s}$ for some $c_s \in\Q\etoile$. Since $\pi$ is transcendental, so is $\zeta(s) $ in this case. No such formula is known, or even conjectured to exist, when $s\geq 3$ is odd. Eventhough $\pi$, $\zeta(3)$, $\zeta(5)$, \ldots are conjectured to be algebraically independent over $\Q$, very few results are known in this direction.
 
 The first one is due to Ap\'ery~\cite{Apery}: $\zeta(3)$ is irrational. Then the next breakthrough is the following result of Ball-Rivoal~\cite{BR, RivoalCRAS}:
 \begin{equation} \label{eqBR}
\dim_\Q \Span_\Q (1, \zeta(3), \zeta(5), \ldots, \zeta(s)) \geq \frac{1-\eps }{1+\log 2}\log s 
 \end{equation}
for any $\eps>0$, provided that $s$ is an odd integer large enough in terms of $\eps$. This result has been made effective, and refined, by several authors -- but only for small values of $s$, and there is still no odd $s\geq 5$ for which $\zeta(s)$ is known to be irrational. For large values of $s$, the following result is the first improvement\footnote{After this paper was written, Lai~\cite{Lai2024} refined the constant $\frac{1 }{1+\log 2}=0.59\ldots$ in Eq.~\eqref{eqBR} to $0.66\ldots$} on the lower bound~\eqref{eqBR}.
\begin{Th} \label{thintroun} For any sufficiently large odd integer $s$ we have:
$$\dim_{\Q } \Span_{\Q } (1,\zeta(3),\zeta(5),\ldots,\zeta(s)) \geq 0.21 \frac{\sqrt{s}}{\sqrt{ \log s}} .$$
\end{Th}
Here $0.21$ is the rounded value of a real number that we did not try to compute exactly.

\bigskip

As a corollary, there are at least $ 0.21 \frac{\sqrt{s}}{\sqrt{ \log s}} $ irrational numbers among $\zeta(3)$, $\zeta(5)$, \ldots, $\zeta(s)$. This weaker result was proved recently by Lai and Yu~\cite{LaiYu} with a better numerical constant, namely\footnote{This constant $1.19\ldots$ has also be refined by Lai~\cite{Lai2025} to $1.28\ldots$} $1.19\ldots$ instead of $0.21$, by following the approach of~\cite{Zudilintrick} and~\cite{Sprang}, developed in~\cite{FSZ}. This strategy provides only a lower bound on the number of irrational odd zeta values, but nothing like~\eqref{eqBR} or Theorem~\ref{thintroun} about linear independence. This makes an important difference: no linear independence criterion is needed, so that the proof is much more elementary.

\bigskip

The proof of Theorem~\ref{thintroun} extends to values of polylogarithms $\Li_s(z) = \sum_{n=1}^\infty \frac{z^n}{n^s}$; recall that $\Li_1(z) = -\log(1-z)$. From now on, we fix an embedding of $\Qbar $ in $\C$. 
Given a positive integer $s$, and $z\in\Qbar\etoile$ such that $|z|$ is small enough (in terms of $s$ and the degree and height of $z$), the values $1$, $\Li_1(z)$, \ldots, $\Li_s(z)$ are known to be $\Q(z)$-linearly independent 
(see~\cite{Nikishin, Hatapolylogs} for the case $z\in\Q$, and~\cite{Chudseul, Chuddeux, Andre} for the general case).
If $z\in\Qbar\etoile$ is fixed with $|z|<1$, this is conjecturally true for any $s$ but the only known result is the following one (due to Rivoal~\cite{Rivoalpolylogs} for $z\in\R$, to Marcovecchio~\cite{Marcovecchio} in the general case): for any non-zero $z\in\Qbar$ such that $|z|<1$ we have 
$$\dim_{\Q(z)} \Span_{\Q(z)} (1,\Li_1(z),\ldots,\Li_s(z))\geq \frac{1 -\eps}{(1+\log 2)[\Q(z):\Q]} \log s $$
provided $s\in\N=\{0,1,2,\ldots\}$ is sufficiently large in terms of $\eps>0$. We refer also to~\cite{gfndio2} for algebraic points $z$ outside the unit disk.

In this paper we improve this lower bound as follows.

\begin{Th} \label{thpolylogs} Let $s $ be a sufficiently large integer. 
Then for any $z\in\Qbar$ such that $|z|\leq1$ and $z\not\in\{0,1\}$ we have:
$$\dim_{\Q(z)} \Span_{\Q(z)} (1,\Li_1(z),\Li_2(z),\ldots,\Li_s(z))\geq \frac{0.26}{[\Q(z):\Q]} \frac{\sqrt{s}}{\sqrt{ \log s}} .$$
\end{Th}

Of course this result holds trivially at $z=1$ (after removing $\Li_1(z)$ from the family), since even powers of $\pi$ are linearly independent over $\Q$.

\bigskip

Most proofs of irrationality (or linear independence) of odd zeta values start with a rational function 
$$
F_n(X) = \sum_{i=1}^a \sum_{j=0}^n \frac{\cij}{(X+j)^i} \in\Q(X)
$$
where $\cij\in\Z$. For instance Ball-Rivoal's proof of~\eqref{eqBR} is based on the following function
(where $n$ is even and $s$ is odd), which is related to a well-poised hypergeometric series:
$$
F_n(X) = d_n^s n!^{s-2r} \frac{ (X-rn)_{rn} (X+n+1)_{rn}}{(X)_{n+1}^s},
$$
where $(x)_\alpha = x(x+1)\ldots (x+\alpha-1)$ is Pochhammer's symbol, $d_n = \lcm(1,2,\ldots, n)$, and $r = \lfloor \frac{s}{(\log s)^2}\rfloor$. The point is to obtain a linear combination of $1$ and odd zeta values, namely 
\begin{equation} \label{eqflintro}
\sum_{t=1}^\infty F_n(t) = \varrho_{0,n} + \varrho_{3,n}\zeta(3) + \varrho_{5,n}\zeta(5) \ldots+ \varrho_{s,n}\zeta(s)
\eneq
with $ \varrho_{i,n} \in\Z$, such that $| \varrho_{i,n}| \leq \beta^{n(1+o(1))}$ as $n\to\infty$ and the absolute value of~\eqref{eqflintro} is less than $\alpha^{n(1+o(1))}$. Applying a linear independence criterion yields a lower bound $1-\frac{\log \alpha}{\log\beta}$ on the dimension of the $\Q$-vector space spanned by $1$, $\zeta(3)$, $\zeta(5)$, \ldots, $\zeta(s)$.

\bigskip

In the literature, this strategy has always been applied to an explicit rational function $F_n(X)$ with explicit integers $\cij$. This has allowed Ball-Rivoal to bound from below the absolue value of~\eqref{eqflintro}, and apply Nesterenko's linear independence criterion~\cite{Nesterenkocritere}. 

On the contrary, to prove Theorem~\ref{thintroun} we apply Siegel's lemma and obtain in this way the existence of integers $\cij$, not all zero, satisfying suitable assumptions. These integers are therefore {\em not explicit}. This allows us to get completely different asymptotic values of the parameters as $s\to\infty$. Whereas $\log\alpha\sim -s\log s$ and $\log\beta\sim (1+\log 2)s$ in Ball-Rivoal's proof, we obtain $\log\alpha\sim - 4.55 \sqrt{s \log s} $ and $ \log\beta\sim 20.93 \log s$. In particular the coefficients $\cij$ are much smaller than in explicit constructions. 

Using non-explicit integers $\cij$ makes it impossible to use Nesterenko's linear independence criterion. We use Siegel's criterion instead, by considering for each $n$ a family of linear forms instead of just~\eqref{eqflintro}. This extrapolation procedure is performed using derivation with respect to both $t$ and $z$ (see parameters $p$ and $k$ in \S \ref{subsecdio1}). Then a multiplicity estimate (namely a generalization~\cite{SFcaract} of Shidlovsky's lemma) is used to provide sufficiently many linearly independent linear forms. Since $z=1$ is a singularity of the underlying differential system, we work at the point $z=-1$ by taking profit of the classical relation $\Li_i(-1) = (2^{1-i}-1)\zeta(i)$ for $i\geq 2$. In such a setting, for each $n$ multiplicity  estimates usually give $p$ linearly independent linear forms in $p$ numbers. However, in our situation it is not always possible to obtain this: the conclusion of our multiplicity  estimate is weaker, but sufficient because we use a refinement of  Siegel's  linear independence criterion.

 \bigskip
 
The structure of this paper is as follows. Section~\ref{secoutils} contains the tools we need: a version of Siegel's lemma combining equalities and inequalities, a refined version of Siegel's  linear independence criterion, and a generalization of Shidlovsky's lemma. 
In \S \ref{sectech} we apply Siegel's lemma to construct the integers $\cij$, or in other words the rational function $F_n(X)$, that will allow us to prove 
Theorems~\ref{thintroun} and~\ref{thpolylogs} in \S \ref{secdio}.

\section{Diophantine tools} \label{secoutils}

We gather in this section the auxiliary Diophantine tools we shall use in the proof of Theorems~\ref{thintroun} and~\ref{thpolylogs}, namely Siegel's lemma, a refined version of Siegel's linear independence criterion, and a multiplicity estimate which generalizes Shidlovsky's lemma.

\subsection{Siegel's lemma} \label{subseclemLSenonce}

We shall apply the following version of Siegel's lemma. The difference with respect to usual statements (see for instance~\cite[Chapter 1, Lemmas 1, 4D or 9A]{SchmidtLNM}) is that linear inequalities (namely~\eqref{eqlemLS2} below) appear: there are not only linear equations with integer coefficients.

\begin{Lem} \label{lemLS} 
Let $N>M\geq M_0\geq 0$ be integers, and $\lambda_{i,m}\in\Z$ for $1\leq i \leq N$ and $1\leq m \leq M$. For each $1\leq m\leq M$, let $H_m\geq 1 $ be a real number such that 
 $\sqrt{\sum_{i=1}^N \lambda_{i,m} ^2}\leq H_m$. For each $m$ such that $M_0< m\leq M$, let $G_m\geq 1 $ be a real number. Define
 $$X = \sqrt N \Big( H_1\ldots H_{M_0}G_{M_0+1}\ldots G_M\Big)^{\frac1{N-M_0}}.$$
Then there exists $(x_1,\ldots,x_N)\in\Z^N\setminus\{(0,\ldots,0)\}$ such that 
\begin{equation} \label{eqlemLS1}
 \sum_{i=1}^N \lambda_{i,m} x_i =0 \mbox{ for any } m\in\{1,\ldots,M_0\},
\eneq
\begin{equation} \label{eqlemLS2}
\Big| \sum_{i=1}^N \lambda_{i,m} x_i \Big| \leq \frac{H_m X}{G_m} \mbox{ for any } m\in\{M_0+1,\ldots,M\},
\eneq
and
\begin{equation} \label{eqlemLS3}
 \sqrt{\sum_{i=1}^N x_i^2}\, \, \leq \, X.
\eneq
\end{Lem}

Inequality~\eqref{eqlemLS2} means that the upper bound deduced from~\eqref{eqlemLS3} using Cauchy-Schwarz inequality 
 is improved by a multiplicative factor $1/G_m$.

In applying Lemma~\ref{lemLS} we shall use the following consequence of~\eqref{eqlemLS3}:
$$
|x_i| \leq X \mbox{ for any } i\in\{1,\ldots,N\}.
$$

\bigskip

\Dem of Lemma~\ref{lemLS}:
Let $F$ denote the set of all $x = (x_1,\ldots,x_N)\in\R^N$ such that~\eqref{eqlemLS1} holds: this is a Euclidean space of dimension $D\geq N-M_0$, with norm given by $\Vert x\Vert = \sqrt{\sum_{i=1}^N x_i^2} $. It is rational, i.e. given by linear equations~\eqref{eqlemLS1} with integer coefficients $\lambda_{i,m}$; this is equivalent to the existence of a basis of $F$ consisting in elements of $\Q^N$. Then $\Lambda = F \cap \Z^N$ is a lattice in $F$, that is a discrete $\Z$-module of rank $D$; we refer to~\cite[Chapter 1]{SchmidtLNM} for all notions of geometry of numbers used in this proof. We point out that geometry of numbers is considered, in~\cite{SchmidtLNM} and in most references, in the Euclidean space $\R^D$. Since we need to work in $F$, which is Euclidean with the scalar product induced from the canonical one on $\R^N$, we fix a linear isometric isomorphism $F\to \R^D$ and use it to carry all definitions and properties.

The determinant of $\Lambda$, denoted by $\det \Lambda$, is the absolute value of the determinant of any $\Z$-basis of $\Lambda$ with respect to an orthonormal basis of $F$ (because such an orthonormal basis is mapped to the canonical basis of $\R^D$ by the above-mentioned isometric isomorphism). It is equal to the volume of the fundamental parallelepiped of $\Lambda$ (see~\cite[Chapter 1, \S 2]{SchmidtLNM}). 

The {\em height} of $F$, denoted by $H(F)$, is by definition $\det \Lambda$ (see~\cite[Chapter 1, \S 4]{SchmidtLNM} or~\cite{SchmidtAnnals}). Now let $F^\perp$ denote the orthogonal complement of $F$ in $\R^N$, and consider the vector $u_m=(\lambda_{1,m}, \ldots, \lambda_{N,m})\in\Z^N$ for any $m\in\{1,\ldots,M_0\}$. The definition~\eqref{eqlemLS1} of $F$ implies $F^\perp = \Span(u_1,\ldots,u_{M_0})$. Reindexing $u_1$, \ldots, $u_{M_0}$ if necessary, we may assume that $u_1$, \ldots, $u_{N-D}$ are linearly independent, so that $F^\perp = \Span(u_1,\ldots,u_{N-D})$. Denoting by $U$ the square matrix of size $N-D$ of which the columns are the coordinates of $u_1$, \ldots, $u_{N-D}$ in an orthonormal basis of $F^\perp$, since $F^\perp\cap \Z^N$ contains the $\Z$-module spanned by $u_1$, \ldots, $u_{N-D}$ we have 
$$ H(F^\perp) = \det( F^\perp\cap \Z^N)\leq | \det U | \leq \prod_{m=1}^{N-D} \Vert u_m\Vert \leq \prod_{m=1}^{N-D} H_m$$
using Hadamard's inequality (as in~\cite[Chapter 1, \S 4, p. 11]{SchmidtLNM}). Since $H(F)=H(F^\perp)$ (see~\cite[Lemma 4C]{SchmidtLNM}) and $H_m\geq 1$ for any $m$, we have
\begin{equation} \label{eqref1}
\det \Lambda = H(F) \leq \prod_{m=1}^{M_0} H_m.
\eneq

Now let us denote by $\calC$ the set of all $x = (x_1,\ldots,x_N)\in F$ such that Eqns.~\eqref{eqlemLS2} and~\eqref{eqlemLS3} hold. We claim that
\begin{equation} \label{eqref2}
\vol \calC \geq \frac{(2X/\sqrt D)^D}{\prod_{m=M_0+1}^{M} G_m} 
\eneq
where $\vol \calC$ is the volume of $\calC$ inside the Euclidean space $F$. Admitting this lower bound for now, and comparing it with Eq.~\eqref{eqref1} and the definition of $X$, we obtain
$$
\vol \calC \geq 2^D \prod_{m=1}^{M_0} H_m\geq 2^D\det \Lambda
$$
since $N-M_0\leq D \leq N$ and $H_m, G_m\geq 1$ for any $m$. Now $\calC$ is a symmetric compact convex body, so Minkowski's first theorem asserts the existence of a non-zero $x\in\calC\cap\Lambda=\calC\cap\Z^N$. This concludes the proof of Lemma~\ref{lemLS}, except for the claim~\eqref{eqref2} that we shall prove now.

To prove Eq.~\eqref{eqref2} we consider $u_m= (\lambda_{1,m}, \ldots, \lambda_{N,m})$ for any $m\in\{M_0+1,\ldots,M \}$, and notice that $\calC$ contains all $x\in F$ such that 
$$ \Vert x \Vert \leq X \mbox{ and }|\langle u_m, x \rangle | \leq \frac{ \Vert u_m \Vert \, X}{G_m}\mbox{ for any } m\in\{M_0+1,\ldots,M \}$$
since $ \Vert u_m \Vert \leq H_m$. Now all indices $ m\in\{M_0+1,\ldots,M \}$ play symmetric roles so we may assume that $G_{M_0+1}\geq \ldots\geq G_M\geq 1$. There exists an orthonormal basis $(e_1,\ldots,e_D)$ of $F$ such that $u_{M_0+i}\in \Span(e_1,\ldots,e_i)$ for any $1\leq i \leq M-M_0$. We shall prove that $\calC$ contains the set $\calC'$ of all points $x=\alpha_1e_1+\ldots+\alpha_De_D$ such that 
$$ |\alpha_i| \leq \frac{ X}{G_{M_0+i} \sqrt{D}} \quad \mbox{ if } 1\leq i \leq M-M_0, \quad \mbox{ and } \quad |\alpha_i| \leq \frac{ X}{\sqrt{D}} \quad \mbox{ if } M-M_0+1\leq i \leq D.$$
Indeed any such $x$ satisfies $\Vert x \Vert \leq \sqrt D \max_{1\leq i \leq D} |\alpha_i| \leq X$. Moreover, for $M_0+1\leq m \leq M$ we have, since $u_{m}\in \Span(e_1,\ldots,e_{m-M_0})$:
$$  | \langle u_m, x \rangle | =  | \langle u_m, \sum_{i=1}^{m-M_0} \alpha_i e_i \rangle | \leq \Vert u_m \Vert \cdot \Vert \sum_{i=1}^{m-M_0} \alpha_i e_i \Vert \leq \Vert u_m \Vert \sqrt{m-M_0} \frac{ X}{G_m\sqrt{D}} \leq \frac{ \Vert u_m \Vert \, X}{ G_m}.$$
Thus $\calC' \subset \calC $, and Eq.~\eqref{eqref2} follows. 
This concludes the proof of Lemma~\ref{lemLS}.

\subsection{A refinement of Siegel's linear independence criterion} \label{subsecsiegel}

The proof of Theorems~\ref{thintroun} and~\ref{thpolylogs} relies on the following refinement of Siegel's linear independence criterion (for usual versions, see for instance~\cite[p. 81--82 and 215--216]{EMS},~\cite[\S 3]{Matala-Aho},~\cite[Proposition 4.1]{Marcovecchio},~\cite[Proposition 4.6]{SFcaract} 
or~\cite[Theorem 4]{gfndio2}). 

Let $\K$ be a number field embedded in $\C$, and $\OK$ be its ring of integers. Let $\K_\infty=\R$ if $\K\subset\R$, and $\K_\infty=\C$ otherwise. The house of $\xi\in\K$, denoted by $\house{\xi}$, is the maximum modulus of the conjugates of $\xi$.

\begin{Prop} \label{propsiegel}
Let $\theta_0,\ldots,\theta_p$ be elements of $\K_\infty$, with $\theta_0\neq 0$. Let $\tau>0$,  and $(Q_n)$ be a sequence of real numbers with limit $+\infty$. Let $\calN$ be an infinite subset of $\N$, and for any $n\in\calN$ let $ [\ell_{i,j}^{(n)}]_{0\leq i \leq I_n, 0\leq j \leq p}$ be a matrix with coefficients in $\OK$  such that:
\begin{itemize}
\item[$(i)$]
As $n\to\infty$ with $n\in\calN$, 
$$\max_{  i,j  } \house{ \ell_{i,j}^{(n)} } \leq Q_n^{1+o(1)}\quad 
\mbox{ and } \quad 
\max_{0\leq i\leq I_n } | \ell_{i,0}^{(n)} \theta_0 + \ldots + \ell_{i,p}^{(n)} \theta_p | \leq Q_n^{-\tau + o(1)}.$$
\item[$(ii)$] For any $n\in\calN$ sufficiently large, for any $x_0,\ldots,x_p\in\K$, if 
$$\forall i\in\{0,\ldots,I_n\} \quad\quad   \ell_{i,0}^{(n)} x_0 + \ldots + \ell_{i,p}^{(n)} x_p = 0$$
then $x_0=0$.
\end{itemize}Then we have
$$\dim_\K\Span_\K(\theta_0,\ldots,\theta_p)\geq \frac{ [\K_\infty : \R]}{ [\K:\Q]} \cdot ( \tau+1) .$$
\end{Prop}

Usually, in  Siegel's linear independence criterion the conclusion of assumption $(ii)$ is $x_0=\ldots=x_p=0$. It turns out that $x_0=0$ is sufficient (since we assume $\theta_0\neq 0$). In terms of the matrix  $ [\ell_{i,j}^{(n)}]_{0\leq i \leq I_n, 0\leq j \leq p}$, we assume that the first column  (corresponding to $j=0$) is not a linear combination of other  columns. This is weaker than the usual assumption, namely that the $p+1$  columns are linearly independent. In particular, Proposition \ref{propsiegel} may apply in settings where the number of linear forms, namely $I_n+1$, is less than $p+1$.

Another way of stating this is that we assume that the $I_n+1$ linear forms we have (for a given $n$) have no common zero $x=(x_0,\ldots,x_p)$ such that $x_0\neq 0$. The usual assumption is that they have no common zero at all in $\K^{p+1}\setminus\{0\}$. 

\bigskip

In the proof of Theorem~\ref{thintroun} we apply Proposition \ref{propsiegel} with $\K=\Q $, $Q_n = \beta^n$, and $\tau = -\frac{\log \alpha}{\log\beta}$ (so that $Q_n^{-\tau} = \alpha^n$), where $\alpha $ and $\beta$ will be defined in \S \ref{subsecfin}. The setting is similar for 
Theorem~\ref{thpolylogs}, with $\K=\Q(z)$ (see \S \ref{subsecpolylogs}).

\begin{proof}[Proof of Proposition \ref{propsiegel}]
Let $d=\dim_\K \Span_\K(\theta_0,\ldots,\theta_p)$. There exists a matrix \linebreak $\Lambda=[\lambda_{i,j}]_{1\leq i \leq p+1-d, 0\leq j \leq p}\in M_{p+1-d,p+1}(\OK)$ of rank $p+1-d$ such that for any $1\leq i \leq p+1-d$, we have $\sum_{j=0}^{p} \lambda_{i,j}\theta_j=0$. For any $n\in \calN $ we let $L^{(n)}=[\ell_{i,j}^{(n)}]_{0\leq i \leq I_n, 0\leq j\leq p}\in M_{I+1,p+1}(\OK)$ and consider the matrix $M^{(n)}=\left[\begin{array}{c} \Lambda \\  L^{(n)}\end{array}\right]\in M_{p+I_n+2-d, p+1}(\OK)$. Since $\rk M^{(n)}$ takes only finitely many values, there exists $r\leq p+1$ such that $ \rk M^{(n)} = r$ for infinitely many $n\in\calN$. Discarding other elements of $\calN$, we may assume that $ \rk M^{(n)} = r$ for any $n\in\calN$. 

For any $n\in\calN$ there exists $\calC_n\subset\{0,\ldots,p\}$ of cardinality $r$ such that the columns of $M^{(n)} $ with index $j\in \calC_n$ are linearly independent; then all other columns of $M^{(n)} $  are $\K$-linear combinations of these $r$ columns. If $0\not\in\calC_n$ then there exist $x_1^{(n)}, \ldots, 
x_p^{(n)}\in\K$ such that for any $0\leq i \leq I_n$, $\ell_{i,0}^{(n)} = \sum_{j=1}^p x_j^{(n)}\ell_{i,j}^{(n)}$. This contradicts hypothesis $(ii)$ if $n$ is large enough, so that $0\in\calC_n$ (discarding finitely many integers $n$ if necessary).  Since $\calC_n$  takes only finitely many values, as above we may assume that there exists  $\calC\subset\{0,\ldots,p\}$ of cardinality $r$, with $0\in\calC$, such that $\calC_n=\calC$ for any $n\in\calN$. Since $\theta_1,\ldots,\theta_p$ play symmetric roles, we assume for simplicity that $\calC = \{0,\ldots, r-1\}$.

We denote by $C_0^{(n)},\ldots, C_p^{(n)}$ the columns of $M^{(n)}$. Since $\calC_n = \{0,\ldots, r-1\}$,  for any $j\in \{r,\ldots,p\}$ there exist $\kappa_{j,0}^{(n)}, \ldots, \kappa_{j,r-1}^{(n)}\in\K$ such that $C_j^{(n)} = \kappa_{j,0}^{(n)} C_0^{(n)} + \ldots+  \kappa_{j,r-1}^{(n)} C_{r-1}^{(n)} $. This implies $ \kappa_{j,0}^{(n)} \ell_{i,0}^{(n)} + \ldots+  \kappa_{j,r-1}^{(n)}  \ell_{i,r-1}^{(n)} -  \ell_{i,j}^{(n)}=0$ for any $0\leq i \leq I_n$, so that $ \kappa_{j,0}^{(n)}=0$ using assumption $(ii)$. We deduce that for any $0\leq i \leq I_n$:
$$
\sum_{j=0}^p  \ell_{i,j}^{(n)}\theta_j= \sum_{j=0}^{r-1}  \ell_{i,j}^{(n)}\theta_j +  \sum_{j=r}^{p}\Big( \sum_{t=1}^{r-1}   \kappa_{j,t}^{(n)}  \ell_{i,t}^{(n)} \Big)  \theta_j  =  \sum_{s=0}^{r-1} \Theta_s^{(n)}  \ell_{i,s}^{(n)}
$$
where 
$$
\Theta_s^{(n)} = \theta_s +  \sum_{j=r}^{p}    \kappa_{j,s}^{(n)}  \theta_j  \quad \mbox{ for any $0\leq s \leq r-1$,}
$$
and in particular $\Theta_0^{(n)} = \theta_0$ since  $ \kappa_{j,0}^{(n)}=0$ for any $j$. In the same way, we have $0 = \sum_{j=0}^p \lambda_{i,j}\theta_j =  \sum_{s=0}^{r-1} \Theta_s^{(n)}  \lambda_{i,s} $ for any $1\leq i \leq p+1-d$. Therefore the linear combination of columns $\sum_{s=0}^{r-1} \Theta_s^{(n)}  C_s^{(n)}$ consists in $p+1-d$ coefficients equal to $0$, and then $I_n+1$ coefficients bounded by $Q_n^{-\tau + o(1)}$. 

Denote by $M_1^{(n)}=\left[\begin{array}{c} \Lambda_1 \\  L_1^{(n)}\end{array}\right]\in M_{p+I_n+2-d, r}(\OK)$ the matrix obtained by keeping only the first $r$ columns of  $M^{(n)}=\left[\begin{array}{c} \Lambda \\  L^{(n)}\end{array}\right]$, so that $\Lambda_1 $ and $L_1^{(n)}$ are obtained in the same way from  $\Lambda  $ and $L ^{(n)}$ respectively. Then $\rk M_1^{(n)} = r$ since $\calC_n  = \{0,\ldots, r-1\}$, and $\rk \Lambda_1  = \rk \Lambda  = p+1-d$ because the columns of $\Lambda_1$ span the same vector space as those of $\Lambda$; in particular, $r\geq p+1-d$.

Therefore  $M_1^{(n)}=\left[\begin{array}{c} \Lambda_1 \\  L_1^{(n)}\end{array}\right]$ has rank $r$ equal to its number of columns, and  its first $p+1-d$ rows are linearly independent:  we may choose $r-(p+1-d)$ rows among those of  $L_1^{(n)}$ that make up, together with $ \Lambda_1$, an invertible matrix. Up to renumbering the linear forms, we may assume that the first $r-(p+1-d)$ rows have this property. Then we denote by $L_2^{(n)}=[\ell_{i,j}^{(n)}]_{0\leq i \leq  r-p+d-2, 0\leq j\leq r-1}\in M_{r-p+d-1,r}(\OK)$ the matrix obtained from $L_1^{(n)}$ by keeping only these rows, and we let $M_2^{(n)}=\left[\begin{array}{c} \Lambda_1 \\  L_2^{(n)}\end{array}\right]\in M_{r}(\OK)\cap {\rm GL}_r(\K)$.

As in the usual proof of Siegel's criterion, we may now consider the non-zero determinant $\Delta^{(n)}\in \OK$ of $M_2^{(n)}$. Recall that $\K$ is  embedded in $\C$, and that  $\sum_{s=0}^{r-1} \Theta_s^{(n)}  C_s^{(n)}$ consists in $p+1-d$ coefficients equal to $0$, and then $I_n+1$ coefficients bounded by $Q_n^{-\tau + o(1)}$, where $C_0^{(n)},\ldots, C_{r-1}^{(n)}$ are the columns of $M_1^{(n)}$. Keeping only the first $r$ coefficients of these columns, we obtain the corresponding columns of $M_2^{(n)}$. Then  $\Delta^{(n)}$ is equal to the determinant of the matrix obtained from  $M_2^{(n)}$ by replacing   the first column   with this linear combination divided by $\Theta_0^{(n)} = \theta_0$  (which is non-zero by assumption, and independent of $n$). Since only the last  $r-p+d-1 $ rows of $M_2^{(n)}$ depend on $n$, and these are  the only rows where non-zero coefficients may appear in the new first column, we obtain by expanding the determinant with respect to this column:
$$|\Delta^{(n)}| \leq  Q_n^{-\tau + o(1)} \Big(Q_n^{1+o(1)}\Big)^{r-p-2+d} =  Q_n^{-\tau +r-p-2+d+ o(1)} \leq  Q_n^{-\tau  -1+d+ o(1)} 
$$ 
using assumption $(i)$ and the upper bound $r\leq p+1$.

Let $\delta=[\K:\Q]$ and denote by $\sigma_1=\Id$, $\sigma_2$, \ldots, $\sigma_\delta$ the embeddings $\K\to\C$. If $\K_\infty=\R$, we bound $|\sigma_k(\Delta^{(n)})| = |\det \sigma_k(M_2^{(n)})|$ trivially by $Q_n^{r-p-1+d+o(1)}\leq Q_n^{d+o(1)} $ for any $2\leq k \leq \delta$, so that $ \prod_{k=1}^\delta  \sigma_k(\Delta^{(n)})  \in\Z\setminus\{0\}$ satisfies
$$1 \leq \Big| \prod_{k=1}^\delta \sigma_k(\Delta^{(n)}) \Big| \leq Q_n^{-\tau  -1+d+d(\delta-1) +o(1)} = Q_n^{-\tau  -1+d \delta + o(1)}$$
and therefore $d \delta \geq \tau+1$. If  $\K_\infty=\C$ then we may assume $\sigma_2$ to be complex conjugation so that $ |\sigma_2(\Delta^{(n)})| =  |\sigma_1(\Delta^{(n)})|$. We bound  $|\sigma_k(\Delta^{(n)})| $ as above for $3\leq k \leq \delta$, and deduce  $2( -\tau  -1+d) + (\delta-2)d\geq 0$, that is $d \delta \geq 2(\tau+1)$.
\end{proof}

\subsection{Multiplicity estimate} \label{subsecshidenonce}

Let us state now the generalisation of Shidlovsky's lemma we shall use, namely~\cite[Theorem~3.1]{SFcaract}. It is based on Fuchs' global relation on exponents, following the approach initiated by Chudnovsky~\cite{ChudShid,ChudShidDeux} in the Fuchsian case and generalized by 
Bertrand-Beukers~\cite{BB} and Bertand~\cite{DBShid} using differential Galois theory.

\bigskip

We consider a positive integer $\qshid$ and a matrix $A \in M_\qshid(\C(z))$. We let $S_0,\ldots,S_{\qshid-1}\in\C[X]$ with $\deg S_i \leq m$ for any $i$. With each solution $Y = \tra (y_0,\ldots,y_{\qshid-1})$ of the differential system $Y'=AY$ is associated a remainder $R(Y)$ defined by 
$$R(Y)(z) = \sum_{i=0}^{\qshid-1} S_i(z) y_i(z).$$
Let $\Sigma$ be a finite subset of $\P^1(\C) = \C\cup\{\infty\}$, with $\infty\in\Sigma$. For each $\sigma\in\Sigma$, let 
$(Y_j)_{j\in J_\sigma}$ be a family of solutions of $Y'=AY$ such that:
\begin{itemize}
\item For any $j\in J_\sigma$, the function $R(Y_j)$ belongs to the Nilsson class at $\sigma$, i.e. can be written as a $\C$-linear combination of functions of the form $h(z)(z-\sigma)^a(\log(z-\sigma))^b$ with $a\in\C$, $b\in\N$, and $h$ holomorphic at $\sigma$; here $z-\sigma$ should be understood as $1/z$ if $\sigma=\infty$. 
\item The functions $R(Y_j)$, for $j\in J_\sigma$, are linearly independent over $\C$ (as functions on a small open disk centered at $\sigma$).
\end{itemize}

\begin{Th} \label{thzerofct}
Let $\mu$ denote the order of a non-zero differential operator $L \in \C(z)[\frac{\dd}{\dd z}]$ such that $L(R(Y_j))= 0 $ for any $\sigma\in\Sigma$ and any $j\in J_\sigma$. Then 
$$
\sum_{\sigma\in\Sigma}\sum_{j\in J_{\sigma}} \ord_\sigma(R(Y_j)) \leq (m+1) (\mu - \Card J_\infty) + \cstun
$$
where $\cstun$ is a constant that depends only on $A$ and $\Sigma$.
\end{Th}

In this result we denote by $ \ord_\sigma$ the order of vanishing at $\sigma$ (recall that logarithmic factors may appear, but they have no influence on the order of vanishing; for instance, $\ord_0(z^e(\log z)^i)$ is the real part of $e$, for $e\in\C$ and $i\in\N$). 

\section{A non-explicit rational function} \label{sectech}

In this section we construct the rational function $F_n(X)$ that will be used in \S \ref{secdio} to prove
Theorems~\ref{thintroun} and~\ref{thpolylogs}. The output of this construction is stated as Theorem~\ref{thsiegel} in \S \ref{subsec31}. Its proof, based on Siegel's lemma, is given in \S \ref{subsecdio2}. It relies on a result of~\cite{FR}, which relates asymptotic estimates of $F_n(X)$ at infinity to values at 1 of some functions $P_{k,1}(z)$ related to a differential system arising from polylogarithms. In \S \ref{subsectech1} we define these functions $P_{k,1}(z)$, explain the setting and state as Proposition~\ref{propgene} a technical result used in the proof of Theorem~\ref{thsiegel}. We prove 
 Proposition~\ref{propgene} in \S \ref{subsectech3}, after dealing with a lemma of analytic number theory in \S \ref{subsectech2}.

\subsection{Output of the construction}\label{subsec31}

In this section we 
apply Siegel's lemma (namely Lemma~\ref{lemLS} stated in \S \ref{subseclemLSenonce}) to construct integers $\cij\in\Z$, for $1\leq i \leq a$ and $0\leq j \leq n$, such that the rational function
\begin{equation}\label{eqfncij}
F_n(X) = \sum_{i=1}^a \sum_{j=0}^n \frac{\cij}{(X+j)^i} \in\Q(X)
\eneq
has interesting properties. We denote by 
$$F_n(t) = \sum_{d=1}^\infty \frac{\frakA_d}{t^d}$$
the expansion of $F_n(t)$ as 
 $|t|\to\infty$.

\begin{Th} \label{thsiegel}
Let $a\in\N$ and $\om,\gdom,r\in\Q$ be such that $a>\gdom\geq \om \geq 1$ and $r\geq 1$. Then for any $n\geq 0$ such that $rn, \om n , \gdom n\in\N$
 there exist integers $\cij\in\Z$ for $1\leq i \leq a$ and $0\leq j \leq n$, not all zero, with the following properties:
 \begin{itemize}
 \item[$(i)$] As $|t| \to\infty$, we have $F_n(t) = O(|t|^{-\om n})$.
 \item[$(ii)$] As $n\to \infty$, we have 
 $|\cij| \leq \chi^{n(1+o(1))}$ for any $i$, $j$, with
\begin{equation} \label{eqdefchi} 
\chi = \exp\Big(\frac{ \om \log 2 + 3\om^2 + \om^2 \log (a+1) +\frac12 \gdom^2 \log r }{a-\om }\Big).
\eneq
 \item[$(iii)$] For any $d < \gdom n$ we have $|\frakA_d |\leq r^{d-\gdom n} n^d d^a \chi^{n(1+o(1))}$.
 \end{itemize}
 Moreover in $(ii)$ and $(iii)$ the sequences denoted by $o(1)$ do not depend on $i$, $j$, $d$, and tend to 0 as $n\to\infty$.
\end{Th}

 The upper bound $(iii)$ is interesting only when $\om n \leq d < \gdom n$, since part $(i)$ means 
 $\frakA_d = 0$ for any $d< \om n$. We also point out that, even if it is not explicit in the notation, the integers $\cij$ depend on $a,\om,\gdom, r, n$.

\bigskip

This section is devoted to the proof of Theorem~\ref{thsiegel}; this proof will be completed in \S \ref{subsecdio2}.

\bigskip

A rather easy construction of integers $\cij$ satisfying properties $(i)$ and $(iii)$ of Theorem~\ref{thsiegel} would be to apply Lemma~\ref{lemLS}, translating $(i)$ as $\frakA_d =0$ for any $d<\om n$. However the explicit expression of $\frakA_d $ (see Eq.~\eqref{eqdefad} in \S \ref{subsecdio2}) shows that for $d$ close to $\om n$, the equation $\frakA_d =0$ is of the form $\sum_{i,j} \lambda_{i,j}\cij=0$ with integers $ \lambda_{i,j}$ such that $| \lambda_{i,j}|\leq n^{\om n(1+o(1))}$. Applying Lemma~\ref{lemLS} with such a huge bound would not give as $n\to\infty$ a geometric bound on $|\cij|$ in $(ii)$, and therefore it would not seem possible to derive any Diophantine application. Instead, to prove Theorem~\ref{thsiegel} we translate assertion $(i)$ as $P_{k,1}(1)=0$ for any $k<\om n$ (see \S \ref{subsecdio2}). We shall define these functions $P_{k,1}(z)$ now. 

\subsection{Setting of the proof} \label{subsectech1}

Let $a\geq 1$ and $n\geq 0$. In this section we start with arbitrary real numbers 
 $\cij$, for $1\leq i \leq a$ and $0\leq j \leq n$, which may either be fixed or considered as unknowns.
 We point out that the result of \S \S \ref{subsectech1} to~\ref{subsectech3}, namely Proposition~\ref{propgene} below, will be used 3 times in this paper: in \S \ref{subsecdio2} to prove Theorem~\ref{thsiegel}, in \S \ref{subsecdio3} to prove Lemma~\ref{lemcoeffs}, and in \S \ref{subsecpolylogs} for Theorem~\ref{thpolylogs}. 
 
 We let $P_i(z) = \sum_{j=0}^n \cij z^j $ for $1\leq i \leq a$, and $P_0(z)=0$. We define $P_{k,i}(z)$ for $0\leq i \leq a$ and $k\geq 1$ as follows: $P_{1,i}(z) = P_i(z)$ for any $i$, and for $k\geq 2$:
\begin{equation} \label{eqrecpkijtechnique}
\left\{ \begin{array}{l}
P_{k,i}(z) = P'_{k-1,i}(z) - \frac1{z} P_{k-1,i+1}(z) \mbox{ for } 1\leq i \leq a \\
P_{k,0}(z) = P'_{k-1,0}(z) + \frac{\ate z + \bte}{z(1-z)} P_{k-1,1}(z) 
\end{array}\right.
\eneq
where $P_{k-1,a+1}$ is taken to be the zero function; the motivation for this definition will be given in \S\S \ref{subsecdio2} and~\ref{subsecdio1} (see Eqns.~\eqref{eqpkdef} and~\eqref{eqdefqpk}).
 Here $(\bte,\ate)\in\Z^2$ is fixed; we shall take $(\bte,\ate)=(1,1)$ in the proof of Theorem~\ref{thintroun}, and $(\bte,\ate)=(1,0)$ for Theorem~\ref{thpolylogs}. It is not difficult (as in~\cite[proof of Proposition~4.4]{SFcaract}) to prove that $z^{k-1} P_{k,i}(z)$ is a polynomial of degree at most $n$ for $1\leq i \leq a$, and that $z^{k-1} (1-z)^{k-1} P_{k,0}(z)$ is a polynomial of degree at most $n+k-1$; this follows also from the proof of Proposition~\ref{propgene} below. We define the coefficients $p_{k,i,j}$ by
\begin{equation} \label{eqdefpkij}
\left\{ \begin{array}{l}
 z^{k-1} P_{k,i}(z) = \sum_{j=0}^n p_{k,i,j} z^j \mbox{ if } i\geq 1,\\
z^{k-1} (1-z)^{k-1} P_{k,0}(z)= \sum_{j=0}^{n+k-1} p_{k,0,j} z^j.
\end{array}\right.
\eneq
It is clear that each coefficient $p_{k,i,j} $ is a $\Q$-linear combination of the (fixed or unknown) coefficients $c_{i',j'}$ we have started with to define $P_0$, \ldots, $P_a$. In other words, there exist rational numbers $\theco$ such that for any $k$, $i$, $j$:
\begin{equation} \label{eqch0}
p_{k,i,j} = \sum_{i'=1}^a\sum_{j'=0}^n \theco c_{i',j'}.
\eneq
The point of the next result, which is the main step in the proof of Theorem~\ref{thsiegel}, is to provide a common denominator (depending only on $k$) and an upper bound on these coefficients $\theco$.

\begin{Prop} \label{propgene}
For any $k\geq 1$ there exists a positive integer $\del_k =\del_k (a,n, \bte,\ate)$, which depends only on $k$, $a$, $n$, $\bte$, $\ate$, such that:
\begin{itemize}
\item[$(i)$] We have $\del_k \leq (e^3(a+1))^{ \max(n,k) }$
 provided $n$ is large enough in terms of $a$.
\item[$(ii)$] For any $i$, $j$, $i'$, $j'$ we have $\frac{\del_k}{(k-1)!} \theco \in \Z$. 
\item[$(iii)$] For any $i$, $j$, $i'$, $j'$ we have 
$$\Big| \frac{\del_k}{(k-1)!} \theco \Big| \leq 
\left\{
\begin{array}{l}
k^a 2^n \del_k \mbox{ if } 1\leq i \leq a,\\
\max(| \bte| ,| \ate| )\, k^{a+1} 8^{\max(n,k)}\del_k \mbox{ if } i=0.
\end{array}
\right.
$$
\end{itemize}
\end{Prop}

The first observation is that we have geometric bounds as $n\to\infty$ (with $k<\om n$): this solves the problem raised at the end of \S \ref{subsec31}. Another crucial remark is the dependence with respect to $a$ of the upper bound in $(i)$: it is polynomial in $a$, whereas a direct approach would lead to an exponential bound, thereby ruining the Diophantine application we have in mind. Indeed we recall (see the end of the introduction, or \S \ref{subsecfin} for details) that we plan to construct a linear combination of odd zeta values, with coefficients bounded by $\beta^{n(1+o(1))}$ as $n\to\infty$, where $\beta$ is a polynomial in $a$. To achieve this, the bound in $(i)$ has to be polynomial in $a$. This property 
comes from Lemma~\ref{lemarithgene} below. 

\bigskip

In the proof of Theorem~\ref{thsiegel} we shall not use the case $i=0$ of parts $(ii)$ and $(iii)$, but they will be used in the proof of Lemma~\ref{lemcoeffs} in \S \ref{subsecdio3}.

\subsection{A lemma from analytic number theory} \label{subsectech2}

A crucial step in the proof of Proposition~\ref{propgene} is the use of the following lemma, which is of independent interest.

\begin{Lem} \label{lemarithgene} Let $a, N \geq 1$. Denote by $\Delta_{a,N}$ the least common multiple of all products $N_1\ldots N_{\lebeta}$ where $\lebeta\leq a$ and $N_1$, \ldots, $N_\lebeta$ are pairwise distinct integers between $-N$ and $N$ such that $\max N_i - \min N_i \leq N$. Then as $N\to\infty$ (while $a$ is fixed) we have:
\begin{equation} \label{eqlemarithgene}
 \Delta_{a,N} = \exp \Big( N ( \sum_{j=1}^a \frac1j+o(1))\Big)
 \leq \Big(( a+1) e^{\gamma+o(1)}\Big)^{ N}
\eneq
where $\gamma$ is Euler's constant.
\end{Lem} 

The naive version of this lemma would be to use the upper bound $ \Delta_{a,N} \leq d_N^a$, where $d_N = \lcm(1,2,\ldots,N)$, leading to $ \Delta_{a,N} \leq e^{aN(1+o(1))}$. The dependence in $a$ is much better in Lemma~\ref{lemarithgene} because we use the assumption that $N_1$, \ldots, $N_\lebeta$ are pairwise distinct. 

\bigskip

In the proof we shall use the function $\psi$ defined by $\psi(x) = \sum_{p^e \leq x} \log p$, where the sum is over prime numbers $p$ and positive integers $e$ such that $p^e \leq x$. The least common multiple of $1$, $2$, \ldots, $N$ is then $\exp(\psi(N))$. We recall (see for instance~\cite[Chapter XXII, Theorem 434]{HW}) that
 the prime number theorem yields $\psi(N) = N(1+o(1))$. 

\bigskip

\Dem of Lemma~\ref{lemarithgene}: 
For any prime power $p^e$ we let $f_{a,N}(p^e) = \min( a, \lfloor \frac{N}{p^e} \rfloor)$ and we consider
$$\Delta = \prod_{p ^e \leq N} p^{ f_{a,N}(p^e) }$$
where the product is taken over all pairs $(p,e)$ such that $p$ is a prime number, $e\geq 1$, and $p^e\leq N$. Our goal is to prove that $\Delta_{a,N} =\Delta$. To begin with, we compute for any prime $p\leq N$ the $p$-adic valuation of $\Delta$ as follows:
\begin{equation} \label{eqvpDel}
v_p(\Delta)= \sum_{e=1}^{\lfloor \frac{\log N}{\log p} \rfloor } f_{a,N}(p^e) =a \Big\lfloor \frac{\log (N/a) }{\log p} \Big\rfloor + \sum_{e = \lfloor \frac{\log (N/a)}{\log p} \rfloor +1 }^{\lfloor \frac{\log N}{\log p} \rfloor }\Big\lfloor \frac{N}{p^e} \Big\rfloor .
\end{equation}
Now let us prove that $ \Delta_{a,N} $ divides $\Delta$. Let $p$ be a prime number; we shall prove that $v_p ( N_1\ldots N_{\lebeta} ) \leq v_p( \Delta)$ for any non-zero pairwise distinct integers $N_1$, \ldots, $N_\lebeta$ between $-N$ and $N$, with $\lebeta\leq a$ and $\max N_i - \min N_i \leq N$. Since $|N_i| \leq N$ for each $i$, we have 
\begin{equation} \label{eqvppro}
v_p ( N_1\ldots N_{\lebeta} ) = \sum_{i=1}^\lebeta v_p(N_i) =\sum_{e =1}^{ \lfloor \frac{\log N }{\log p} \rfloor } \Card \Spe
 \end{equation}
where $\Spe =\{i\in\unlebeta, \, v_p(N_i)\geq e\}$. Obviously we have $\Card \Spe\leq \lebeta \leq a$, and 
$$\Card\Spe\leq \Big\lfloor \frac{\max_i N_i - \min_i N_i }{ p^e} \Big\rfloor +1\leq \Big\lfloor\frac{N }{ p^e}\Big\rfloor +1.$$
Moreover if $\Card\Spe = \lfloor\frac{N }{ p^e}\rfloor +1$ then $ \min_i N_i = up^e$ and $ \max_i N_i = vp^e$ with $u,v\in\Z$ such that $v-u = \lfloor\frac{N }{ p^e}\rfloor$. If $u\geq 1$ then $v\geq 1 + \lfloor\frac{N }{ p^e}\rfloor > N/p^e$ so that $vp^e> N$, which is impossible. The same contradiction holds if $v\leq -1$ because in this case $-u \geq 1 + \lfloor\frac{N }{ p^e}\rfloor > N/p^e$. Therefore we have $u\leq 0 \leq v$; since all $N_i$ are non-zero, we obtain $ \Card\Spe\leq \lfloor\frac{N }{ p^e}\rfloor $ and finally $ \Card\Spe\leq f_{a,N}(p^e) $. Combining Eqns.~\eqref{eqvppro} and~\eqref{eqvpDel} concludes the proof that $ \Delta_{a,N} $ divides $\Delta$.

\bigskip

Let us prove now\footnote{For the application we have in mind, an upper bound on $\Delta_{a,N}$ is enough. We provide its exact asymptotics for the sake of completeness.} that $ \Delta $ divides $\Delta_{a,N}$. Let $p$ be a prime number; we shall construct pairwise distinct integers $N_i$ between 1 and $N$ such that $v_p ( N_1\ldots N_a ) =v_p(\Delta)$. We write $e = \lfloor \frac{\log (N/a) }{\log p} \rfloor +1$, so that $p^{e-1}\leq N/a < p^e$, and $k = \lfloor\frac{N }{ p^e}\rfloor $. If $ \lfloor \frac{\log N }{\log p} \rfloor = \lfloor \frac{\log (N/a) }{\log p} \rfloor$ the sum in Eq.~\eqref{eqvpDel} is empty, so that letting $N_i = i p^{e-1}$ for $1\leq i \leq a$ we have 
$v_p ( N_1\ldots N_a ) = a(e-1)=v_p(\Delta)$ since assuming $ \lfloor \frac{\log N }{\log p} \rfloor = \lfloor \frac{\log (N/a) }{\log p} \rfloor$ implies $a<p$ so that $v_p(i)=0$ for any $1\leq i \leq a$. 
Assume now, on the contrary, that $ \lfloor \frac{\log N }{\log p} \rfloor \geq e$. Then we have $p ^e\leq N$ and $k\geq 1$; we let $N_i = i p^e$ for $1\leq i \leq k$, and we pick up $N_{k+1}$, \ldots, $N_a$ among the $\lfloor\frac{N }{ p^{e-1}}\rfloor - \lfloor\frac{N }{ p^e}\rfloor \geq a-k$ integers between $p^{e-1} $ and $N $ with $p$-adic valuation equal to $e-1$. Then for any $i\in\una$ we have $e-1\leq v_p(N_i)\leq \lfloor\frac{\log N }{\log p }\rfloor$, and for any $e'\in\{e, \ldots, \lfloor\frac{\log N }{\log p }\rfloor \}$ the number of indices $i$ such that $
v_p(N_i)\geq e'$ is equal to $ \lfloor\frac{N }{ p^{e'}}\rfloor $. Therefore we have 
$$v_p ( N_1\ldots N_a ) = a(e-1) + \sum_{e'=e}^{\lfloor\frac{\log N }{\log p }\rfloor} \Big\lfloor\frac{N }{ p^{e'}}\Big\rfloor =v_p(\Delta)$$
using Eq.~\eqref{eqvpDel}. Finally, for any prime $p$ we have found pairwise distinct integers $N_i$ between 1 and $N$ such that $v_p(\Delta) = v_p ( N_1\ldots N_a ) $. Therefore $ \Delta $ divides $\Delta_{a,N}$, and equality holds: $ \Delta = \Delta_{a,N}$.

\bigskip

To conclude the proof of Lemma~\ref{lemarithgene}, we use this explicit expression of $ \Delta$ to compute it asymptotically. 
In what follows we denote by $o(1)$ any quantity that tends to 0 as $N\to\infty$, with $a$ fixed. Since $\psi(N) = N(1+o(1))$ as recalled before the proof, we have
\begin{eqnarray*}
\log \Delta 
&=& \sum_{p^e\leq N} f_{a,N}(p^e) \log p\\
&=& \sum_{p^e\leq N/a} a \log p + \sum_{k=1}^{a-1} \sum_{\frac{N}{k+1} < p^e\leq \frac{N}{k} } k \log p \\
&=& a \psi(N/a) + \sum_{k=1}^{a-1} k \Big( \psi(N/k) - \psi(N/(k+1))\Big)\\
&=& a \psi(N/a) + \sum_{k=1}^{a-1} k \psi(N/k) - \sum_{k=2}^{a } (k-1) \psi(N/ k ) \\
&=& a \psi(N/a) + \psi(N) - (a-1) \psi(N/a) + \sum_{k=2}^{a-1} \psi(N/k) \\
&=& \sum_{k=1}^a \psi(N/k) = N\Big( \sum_{k=1}^a 1/k+o(1)\Big).
\end{eqnarray*}
At last, $\sum_{k=1}^a \frac1{k} - \log (a+1) $ is non-decreasing with respect to $a$, and tends to $\gamma $ as $a\to\infty$, so that $\sum_{k=1}^a 1/k \leq \gamma + \log (a+1)$ for any $a$. 
This concludes the proof of Lemma~\ref{lemarithgene}.

\subsection{Proof of Proposition~\ref{propgene}}\label{subsectech3}

In this section we prove Proposition~\ref{propgene} by computing explicitly the coefficients $\theco$. We shall use the following lemma, proved in~\cite{Farhibinomial} using Kummer's theorem on $p$-adic valuations of binomial coefficients. 

\begin{Lem} \label{lembinomial} Let $N$ be a positive integer. The least common multiple of the binomial coefficients $\binom{N}{i}$, $0\leq i \leq N$, is equal to $\frac{d_{N+1}}{N+1}$ where $d_{N+1} = \lcm(1,2,\ldots,N+1)$.
\end{Lem}

We shall use also the following notation.
Given integers $ 0\leq \ell < k$, we denote by $H_{\ell,k}$ the set of all $\hsoul = (h_0,\ldots,h_\ell)\in (\N\etoile)^{\ell+1}$ such that $h_0+\ldots+h_\ell=k$; we let $H_{\ell,k}=\emptyset$ if $\ell\geq k$ or $\ell<0$. In particular we have $H_{0,k} = \{ (k)\}$.

For $\hsoul\in H_{\ell,k}$ and $T\in\Z$, we let 
$$\fctchi(T,k,\hsoul) = \frac{T(T-1)\ldots (T-k+2)}{\prod_{i=0}^{\ell-1} (T+1-\sum_{j=0}^i h_j)}$$
where empty products are taken equal to 1; notice that all factors in the denominator appear also in the numerator, so that $\fctchi(T,k,\hsoul)\in\Z$.
Here and below we agree that if $T = \sum_{j=0}^{i_0} h_j-1$ for some $i_0 \in \{0,\ldots,\ell-1\}$ (which is then unique), then the zero factor $T+1- \sum_{j=0}^{i_0} h_j$ has to be omitted from both products, in the numerator and in the denominator. In precise terms, we then have
$T+2\leq k$ and 
$$
\fctchi(T,k,\hsoul) = (-1)^{k-T}
\frac{T! (k-T-2)!}{ \prod_{0\leq i \leq \ell-1 \atop i\neq i_0} (T+1- \sum_{j=0}^{i} h_j) }.
$$

\bigskip

The proof of Proposition~\ref{propgene} falls into 4 steps.

\bigskip

\noindent {\bf Step 1}: Computation of $\theco$ for $i\geq 1$.

\bigskip

The goal of this step is to prove by induction on $k\geq 1$ that for any $1\leq I \leq a$ and any $0\leq T \leq n$ we have
\begin{equation} \label{eqch1}
\thecoTIT = (-1)^{I-i} \sum_{\hsoul\in H_{I-i,k}} \fctchi(T,k,\hsoul)\quad\quad \mbox{ if } \max(1,I-k+1)\leq i \leq I
\eneq
and $\thecojIT = 0$ otherwise (with $i\geq 1$), namely
\begin{equation} \label{eqch2}
\thecojIT = 0\quad \quad \mbox{ if ($i\geq 1$ and $j\neq T$) or ($i\geq I+1$) or ($1\leq i \leq I-k$).}
\eneq
The value of $\thecoze$, namely with $i=0$, will be computed in Step 2 below.

\bigskip

An equivalent form of Eqns.~\eqref{eqch1} and~\eqref{eqch2} is the following: for any $1\leq i \leq a$ and any $k\geq 1$, we have
\begin{equation} \label{eqch3}
P_{k,i}(z) = \sum_{t=1-k}^{n+1-k} z^t \Big( \sum_{I=i}^{\min(a,i+k-1)} c_{I,t+k-1} (-1)^{I-i} \sum_{\hsoul\in H_{I-i,k}} \fctchi(t+k-1,k,\hsoul) \Big).
\eneq
We shall now prove Eq.~\eqref{eqch3} by induction on $k\geq 1$.

\bigskip

For $k=1$, Eq.~\eqref{eqch3} holds trivially; indeed it reads $P_{1,i}(z) = \sum_{t=0}^{n} c_{i,t}z^t$ since $H_{0,1}=\{(1)\}$ and $\fctchi(t,1,(1))=1$. Let us assume that Eq.~\eqref{eqch3} holds for $k-1$, with $k\geq 2$. We recall that 
$$
P_{k,i}(z) = P'_{k-1,i}(z) - \frac1{z} P_{k-1,i+1}(z) \mbox{ for } 1\leq i \leq a 
$$
with $P_{k-1,a+1}(z)=0$.
Using Eq.~\eqref{eqch3} twice (since it reduces to $0=0$ if $i=a+1$) we obtain:
\begin{eqnarray*}
P_{k,i}(z) 
&=&
\sum_{t=2-k}^{n+2-k} tz^{t-1} \Big( \sum_{I=i}^{\min(a,i+k-2)} c_{I,t+k-2} (-1)^{I-i} \sum_{\hsoul\in H_{I-i,k-1}} \fctchi(t+k-2,k-1,\hsoul) \Big)\\
&& 
- z^{t-1} \Big( \sum_{I=i+1}^{\min(a,i+k-1)} c_{I,t+k-2} (-1)^{I-i-1} \sum_{\hsoul\in H_{I-i-1,k-1}} \fctchi(t+k-2,k-1,\hsoul) \Big).
\end{eqnarray*}
Letting $t'=t-1$ yields
\begin{eqnarray*}
P_{k,i}(z) 
&=&
\sum_{t'=1-k}^{n+1-k} z^{t' } \sum_{I=i}^{\min(a,i+k-1)} c_{I,t'+k-1} (-1)^{I-i} \\
&&
\Big( (t'+1) \sum_{\hsoul\in H_{I-i,k-1}} \fctchi(t'+k-1,k-1,\hsoul)
+ \sum_{\hsoul\in H_{I-i-1,k-1}} \fctchi(t'+k-1,k-1,\hsoul) \Big); 
\end{eqnarray*}
here zero terms have been added (namely $I= i+k-1$ in the first sum, if $i+k-1\leq a$, and $I=i$ in the second term; notice that $H_{k-1,k-1}=H_{-1,k-1}=\emptyset$). To conclude it is enough to check that for any $t$, $I$ such that $1-k\leq t\leq n+1-k$ and $i\leq I \leq \min(a,i+k-1)$ we have
\begin{equation} 
(t+1) \sum_{\hsoul ' \in H_{I-i,k-1}} \fctchi(t +k-1,k-1,\hsoul ' ) + 
 \sum_{\hsoul '' \in H_{I-i-1,k-1}} \fctchi(t +k-1,k-1,\hsoul '' ) \label{eqch4}
\eneq 
$$
= \sum_{\hsoul \in H_{I-i,k}} \fctchi(t +k-1,k ,\hsoul ).
$$
 Indeed let $\hsoul = (h_0,\ldots, h_{I-i})\in H_{I-i,k}$, so that $h_0 + \ldots + h_{I-i} = k$. If $ h_{I-i} \geq 2$ then 
$$
\fctchi(t +k-1,k,\hsoul ) = \frac{(t+k-1)(t+k-2)\ldots (t+1)}{\prod_{\lambda=0}^{I-i-1}(t+k-\sum_{j=0}^\lambda h_j)}
= (t+1)\, \fctchi(t +k-1,k-1,\hsoul ' ) 
$$
where $ \hsoul ' = (h_0,\ldots, h_{I-i-1},h_{I-i} -1 ) \in H_{I-i,k-1}$.
On the other hand, if $h_{I-i}=1$ then for $\lambda= I-i-1$ we have $ t+k-\sum_{j=0}^\lambda h_j = t+1$ so that
$$
\fctchi(t +k-1,k,\hsoul ) = \frac{(t+k-1)(t+k-2)\ldots (t+2)}{\prod_{\lambda=0}^{I-i-2}(t+k-\sum_{j=0}^\lambda h_j)}
= \fctchi(t +k-1,k-1,\hsoul '' ) 
$$
where $ \hsoul '' = (h_0,\ldots, h_{I-i-1} ) \in H_{I-i-1,k-1}$. This concludes the proof of Eq.~\eqref{eqch4}, and by induction that of Eq.~\eqref{eqch3}.

\bigskip

\noindent {\bf Step 2}: Computation of $\theco$ for $i=0$.

\bigskip

In this step we shall prove that for any $k\geq 1$, any $0\leq j \leq n+k-1$, any $1\leq I \leq a$ and any $0\leq T \leq n$ we have
\begin{equation} \label{eqch8}
\vartheta_{k,0,j,I,T} \, \, \, \, = \, \, \, \, \sum_{\eps = 0}^{ 1 } \, \, \ateeps \, \, \, \, 
\sum_{s'=1-k}^{-1} \, \, \, \, \sum_{t'=-s'-k+\eps}^{n-s'-k+\eps} \, \, \, \, (-1)^{j-t'-k+1}
\eneq
$$
\cdot \binom{s'+k-1}{j-t'-k+1} \, \, \sum_{\alpha=-1-s'}^{k-2} \, \, (t'+1)_{s'+\alpha+1} \, \, (s'+\alpha+2)_{-s'-1} \, \, \vartheta_{k-\alpha-1,1,t'+s'-\eps+k,I,T}
 $$
 where the coefficients $ \vartheta_{k-\alpha-1,1,t'+s'-\eps+k,I,T}$ have been computed in Step 1, and $ \ateeps $ comes from Eq.~\eqref{eqrecpkijtechnique}. In Eq.~\eqref{eqch8} and throughout this paper, all binomial coefficients $ \binom{r}{s}$ are considered to be zero if $s<0$ or $s>r$.
 
 With this aim in mind we define functions $\psi_{k,\eps}(z)$ for $k\geq 1$ and $\eps\in\{0,1\}$ by letting $\psi_{1,\eps}(z)=0$ and 
\begin{equation} \label{eqch7}
\psi_{k,\eps}(z) = \psi_{k-1,\eps} ' (z) + z^{\eps-1}(1-z)^{-1} P_{k-1,1}(z)
\eneq
for any $k\geq 2$. Indeed the recurrence relation
$$
P_{k,0}(z) = P'_{k-1,0}(z) + \frac{\ate z + \bte}{z(1-z)} P_{k-1,1}(z) 
$$
with $P_{1,0}(z) = 0$ yields immediately, by induction: 
\begin{equation} \label{eqch9}
P_{k,0}(z) = \sum_{\eps=0}^1 \ateeps \psi_{k ,\eps} (z) \mbox{ for any } k\geq 1.
\eneq
Let us fix $\eps\in\{0,1\}$. Then Eq.~\eqref{eqch7} implies, by induction, 
$$
\psi_{k,\eps}(z) = \sum_{\alpha=0}^{k-2} \Big(\frac{d}{dz}\Big)^{\alpha} \Big( z^{\eps-1}(1-z)^{-1} P_{k-\alpha-1,1}(z)\Big)
$$
for any $k\geq 1$. Recall that 
$$P_{k-\alpha-1,1}(z) = \sum_{t=\alpha+2-k}^{n+\alpha+2-k} p_{k-\alpha-1,1,t+k-\alpha-2}z^t, $$
so that Leibniz' formula yields
$$
\psi_{k,\eps}(z) = \sum_{\alpha=0}^{k-2} \sum_{t=\alpha+2-k}^{n+\alpha+2-k} p_{k-\alpha-1,1,t+k-\alpha-2}
\sum_{\beta=0}^\alpha \binom{\alpha}{\beta}
(t+\eps-\beta)_\beta z^{t+\eps-\beta-1} (\alpha-\beta)! (1-z)^{-1-\alpha+\beta}.
$$
Letting $t'=t+\eps-\beta-1$ and $s'=-1-\alpha+\beta$ we obtain 
$$
\psi_{k,\eps}(z) = \sum_{s' = 1-k}^{-1} \sum_{t' = -s' -k+\eps}^{n-s'-k+\eps} z^{t'}(1-z)^{s'} 
\sum_{\alpha=-1-s' }^{k-2} p_{k-\alpha-1,1,t'+s' + k-\eps}
(t'+1)_{s'+\alpha+1} (s'+\alpha+2)_{-s'-1} .
$$
For $1-k\leq s'\leq -1$ and $ -s' -k+\eps\leq t' \leq n-s'-k+\eps $ we write now
\begin{eqnarray*}
z^{t'} (1-z)^{s'}
&=& (1-z)^{1-k} \sum_{\sigma=0} ^{s'+k-1} (-1)^\sigma z^{\sigma+t'} \binom{s'+k-1}{\sigma}\\
&=&(1-z)^{1-k} \sum_{j=0} ^{n+k-1} (-1)^{j-t'-k+1} z^{ j+1-k } \binom{s'+k-1}{j-t'-k+1}
\end{eqnarray*}
by letting $j = t'+\sigma+k-1$; notice that the values taken by $j$ form actually a subset of $\{0,\ldots,n+k-1\}$, but additional terms are zero because of the above-mentioned convention on binomial coefficients.
Substituting this formula into the expression for $ \psi_{k,\eps}(z) $ and interchanging summations, we obtain
\begin{eqnarray*}
&& \psi_{k,\eps}(z) = 
 (1-z)^{1-k} \sum_{j=0}^{n+k-1} z^{j+1-k} 
\sum_{s' = 1-k}^{-1} \sum_{t' = -s' -k+\eps}^{n-s'-k+\eps} 
(-1)^{j-t'-k+1} \\
&&
\cdot \binom{s'+k-1}{j-t'-k+1} \sum_{\alpha=-1-s' }^{k-2} p_{k-\alpha-1,1,t'+s' + k-\eps}
(t'+1)_{s'+\alpha+1} (s'+\alpha+2)_{-s'-1} .
\end{eqnarray*}
Using Eqns.~\eqref{eqch0} and~\eqref{eqch9} this concludes the proof of Eq.~\eqref{eqch8}.

\bigskip

\noindent {\bf Step 3}: Denominators.

\bigskip

In this step we prove that assertion $(ii)$ of Proposition~\ref{propgene} holds with
$$\del_k = d_k^2 \Delta_{a,\max(k,n)}$$
where $\Delta_{a,\max(k,n)}$ is defined in Lemma~\ref{lemarithgene}. Since $\gamma\leq 1$, the upper bound $(i)$ on $\del_k $ in Proposition~\ref{propgene} follows immediately from Lemma~\ref{lemarithgene} and the prime number theorem (namely, $d_k = \exp(k(1+o(1)))$).

\bigskip

Let us start with the case $i\geq 1$. We shall prove that 
\begin{equation} \label{eqch11}
\frac{d_k \Delta_{a,\max(k,n)}}{(k-1)!} \fctchi(T,k,\hsoul) \in\Z
\eneq
for any $k\geq 1$, $1\leq I \leq a$, $0\leq T \leq n$, $\max(1,I-k+1)\leq i \leq I$ and any $\hsoul = (h_0,\ldots,h_{I-i})\in (\N\etoile)^{I-i+1}$ such that $h_0+\ldots+h_{I-i}=k$.
Using Eq.~\eqref{eqch3} proved in Step 1 and Eq.~\eqref{eqch0}, this is enough to prove assertion $(ii)$ of Proposition~\ref{propgene} for $i\geq 1$ (even in a stronger form, namely with $d_k \Delta_{a,\max(k,n)}$ instead of $\del_k$) .

\bigskip

To prove~\eqref{eqch11}, we recall that 
\begin{equation} \label{eqch10}
\fctchi(T,k,\hsoul) = \frac{T(T-1)\ldots (T-k+2)}{\prod_{\lambda=0}^{I-i-1} (T+1-\sum_{j=0}^\lambda h_j)}. 
\eneq
If $T-k+2> 0$ then
$$
\frac{d_k \Delta_{a,\max(k,n)}}{(k-1)!} \fctchi(T,k,\hsoul) = d_k \binom{T}{k-1} \frac{ \Delta_{a,\max(k,n)}}{\prod_{\lambda=0}^{I-i-1} (T+1-\sum_{j=0}^\lambda h_j)}\in\Z$$
using Lemma~\ref{lemarithgene}, since the $ T+1-\sum_{j=0}^\lambda h_j$ are $I-i\leq a-1$ pairwise distinct integers between $0$ and $T\leq n\leq \max(k,n)$.

If $T-k+2\leq 0$ then a factor vanishes in the numerator of Eq.~\eqref{eqch10}. In proving Eq.~\eqref{eqch11} we may assume that a factor vanishes in the denominator too, namely $ T+1-\sum_{j=0}^{\lambda_0} h_j$, and in this case these factors have to be omitted in Eq.~\eqref{eqch10}; we then have 
\begin{eqnarray*}
&&\frac{d_k \Delta_{a,\max(k,n)}}{(k-1)!} \fctchi(T,k,\hsoul) \\
&=& (-1)^{T-k+2} \frac{d_k}{(k-1)\combi{k-2}{T}} \frac{ \Delta_{a,\max(k,n)}}{ \prod_{0 \leq\lambda\leq I-i-1\atop \lambda \neq \lambda_0 } (T+1-\sum_{j=0}^\lambda h_j)}\in\Z
\end{eqnarray*}
using Lemmas~\ref{lemarithgene} and~\ref{lembinomial}, since the $ T+1-\sum_{j=0}^\lambda h_j$ with $ \lambda \neq \lambda_0$ are $I-i-1\leq a-2$ pairwise distinct integers between $T-k+2\geq -k+2$ and $T\leq n$, with distance at most $k-2$ from one another.

\bigskip

This concludes the proof of assertion $(ii)$ of Proposition~\ref{propgene} for $i\geq 1$; let us study the case $i=0$ now. Using Eq.~\eqref{eqch8} (see Step 2) it is enough to prove that 
$$
\frac{d_k^2 \Delta_{a,\max(k,n)}}{(k-1)!} \, \, \, \, 
(t'+1)_{s'+\alpha+1} \, \, (s'+\alpha+2)_{-s'-1} \, \, \, \, 
 \vartheta_{k-\alpha-1,1,t'+s'-\eps+k,I,T}\, \, \, \, 
\in\, \, \Z
$$
for any $k\geq 1$, $0\leq \eps \leq 1$, $1-k\leq s' \leq -1$, $-s'-k+\eps \leq t' \leq n-s'-k+\eps$, $-1-s'\leq \alpha\leq k-2$, $1\leq I \leq a$ and $0\leq T\leq n$. 
Now we have proved in the first part of Step 3 that for $i\geq 1$, assertion $(ii)$ of Proposition~\ref{propgene} holds with $d_k \Delta_{a,\max(k,n)}$ instead of $\del_k$, so that 
$$
\frac{d_k \Delta_{a,\max(k,n)}}{(k-1-\alpha)!} \vartheta_{k-\alpha-1,1,t'+s'-\eps+k,I,T} \in\Z.
$$
Since we have
$$d_k \frac{ (k-1-\alpha)!}{ (k-1)!} (t'+1)_{s'+\alpha+1} (s'+\alpha+2)_{-s'-1} 
 = \frac{ d_k }{\combi{k-1}\alpha} \combi{ s'+\alpha+1+t'}{t'}\in\Z $$
using Lemma~\ref{lembinomial}, this concludes the proof of assertion $(ii)$ of Proposition~\ref{propgene}.

\bigskip

\noindent {\bf Step 4}: Absolute values.

\bigskip

To conclude the proof of Proposition~\ref{propgene}, let us prove part $(iii)$.
To bound $| \frac{\del_k}{(k-1)!} \thecojIT |$ from above, we begin with the case where $i\geq 1$ and use Eqns.~\eqref{eqch1} and~\eqref{eqch2} proved in Step 1. Whenever $1\leq I \leq a$ and $0\leq T \leq n$ we have $\Card H_{I-i,k}\leq k^{I-i}\leq k^a$ and, for any $\hsoul\in H_{I-i,k}$:
$$\Big| \frac{ \fctchi(T,k,\hsoul) }{(k-1)!} \Big| \leq \binom T {k-1}\leq 2^T\leq 2^n \mbox{ if } T\geq k-1,$$ 
whereas
$$\Big| \frac{ \fctchi(T,k,\hsoul) }{(k-1)!} \Big| \leq \frac{1}{ (k-1) \binom {k-2}T }\leq 1 \mbox{ if } T\leq k-2.$$
Therefore we obtain
\begin{equation}\label{eqmajonvnvun}
 \Big| \frac{\del_k}{(k-1)!} \thecojIT \Big| \leq k^a 2^n \del_k \mbox{ if } i \geq 1.
 \end{equation}

\bigskip

Let us deal now with the case $i=0$, using Eq.~\eqref{eqch8} proved in Step 2. In this sum there are at most $2k(k-1)$ values of the triple $(\eps, s',\alpha)$. For each value, the sum over $t'$ of $\binom{s'+k-1}{j-t'-k+1}$ is bounded by $2^{ s'+k-1}\leq 2^{k-1}$, and we have 
$$
\Big| (t'+1)_{s'+\alpha+1} (s'+\alpha+2)_{-s'-1} \Big| = 
\left\{
\begin{array}{l}
\alpha! \, \binom{t'+s'+\alpha+1}{t'}\leq \alpha! \, 2^n\mbox{ if } t'\geq 0, \\
\\
0 \mbox{ if } t'<0\leq t'+s'+\alpha+1, \\
\\
 \alpha! \, \binom{-t'-1}{ s'+\alpha+1} \leq \alpha! \, 2^{-t'}\leq \alpha ! \, 2^k \mbox{ if } t'+s'+\alpha+1<0.
\end{array}
\right.
$$
Using Eq.~\eqref{eqmajonvnvun} with $k-\alpha-1$ instead of $k$ we deduce that 
\begin{eqnarray*}
\Big| (t'+1)_{s'+\alpha+1} (s'+\alpha+2)_{-s'-1} \frac{1}{(k-1)!} \vartheta_{k-\alpha-1,1,t'+s'-\eps+k,I,T}\Big|
&\leq& \frac{ \alpha! \, (k-\alpha-2)! \, k^a\, 2^{n+\max(n,k)}}{(k-1)!} \\
& \leq & \frac{k^a\,  2^{n+\max(n,k)}}{k-1} 
\end{eqnarray*}
since $\binom{k-2}{\alpha}\geq 1$. 
Therefore Eq.~\eqref{eqch8} yields 
$$
 \Big| \frac{\del_k}{(k-1)!} \thecozeIT \Big| \leq \max( |\bte|, |\ate|)\ k^{a+1} \ 2^{n+k+\max(n,k)} \ \del_k.
$$
This concludes the proof of Proposition~\ref{propgene}.

\subsection{Application of Siegel's lemma} \label{subsecdio2}

In this section we use Proposition~\ref{propgene} to conclude the proof of Theorem~\ref{thsiegel}. The notation is the one of \S \S \ref{subsec31} and~\ref{subsectech1}; the coefficients $c_{i,j}$ are related to the function $F_n(X)$ we are trying to construct by Eq.~\eqref{eqfncij}.

The asymptotic expansion of $F_n(t)$ at infinity reads
\begin{equation} \label{eqfnasy}
F_n(t) = \sum_{d=1}^\infty \frac{\frakA_d}{t^d} \mbox{ for any $t$ such that $|t|>n$,}
\eneq
where the coefficients $\frakA_d$ are given explicitly (see~\cite[Eq. (17)]{FR}) by 
\begin{equation} \label{eqdefad}
\frakA_d = (-1)^d \sum_{i=1}^{\min(a,d)} \sum_{j=0}^n (-1)^i \binom{d-1}{i-1} j^{d-i} \cij \mbox{ for any } d\geq 1.
\eneq
The important point here is that we have also~\cite[Proposition 2]{FR}
\begin{equation} \label{eqrnasy}
R_n(z) = \sum_{d=1}^\infty \frakA_d (-1)^{d-1} \frac{(\log z)^{d-1}}{(d-1)!} \mbox{ for any $z\in\C$ such that $|z-1|<1$}
\eneq
where 
\begin{equation} \label{eqdefrn}
R_n(z) = \sum_{i=1}^a P_i(z) (-1)^{i-1} \frac{(\log z)^{i-1}}{(i-1)!}.
\eneq
As in \S \ref{subsectech1} we consider 
the rational functions $P_{k,i}(z)$ defined by $P_{1,i}(z) = P_i(z)$ and, for any $k\geq 2$,
\begin{equation} \label{eqpkdef}
 P_{k,i}(z) = P'_{k-1,i}(z) - \frac1{z} P_{k-1,i+1}(z) \mbox{ for } 1\leq i \leq a \\
\eneq
where $P_{k-1,a+1}$ is understood as 0; however we are not interested in $P_{k,0}(z)$ here.
Since the derivative of $ (-1)^{i-1} \frac{(\log z)^{i-1}}{(i-1)!} $ is $\frac{-1}{z} (-1)^{i-2} \frac{(\log z)^{i-2}}{(i-2)!} $ if $i\geq 2$, and $0$ if $i=1$, we have 
$$
R_n^{(k-1)}(z) = \sum_{i=1}^a P_{k,i}(z) (-1)^{i-1} \frac{(\log z)^{i-1}}{(i-1)!} \mbox{ for any } k\geq 1 
$$
and in particular
\begin{equation} \label{eqrnderiv}
R_n^{(k-1)}(1) = P_{k,1}(1).
\eneq

\bigskip

Using Eqns.~\eqref{eqfnasy},~\eqref{eqrnasy} and~\eqref{eqrnderiv} we see that the following assertions are equivalent: 
\begin{itemize}
\item[$(i)$] As $|t|\to\infty$, $F_n(t) = O(|t|^{-\om n})$.
\item[$(ii)$] For any $d\in\{1,\ldots,\om n-1\}$, $\frakA_d = 0$.
\item[$(iii)$] As $z\to 1$, $ R_n(z) = O((z-1)^{\om n -1})$.
\item[$(iv)$] For any $k\in \{1,\ldots,\om n-1\}$, $R_n^{(k-1)}(1) = 0$.
\item[$(v)$] For any $k\in \{1,\ldots,\om n-1\}$, $P_{k,1}(1) = 0 $. 
\end{itemize}
Using the notation of \S \ref{subsectech1}, the last assertion reads $\sum_{j=0}^n p_{k,1,j} =0$, or equivalently 
\begin{equation} \label{eqsyslin}
 \frac{\del_k}{(k-1)!} \sum_{i'=1}^a\sum_{j'=0}^n \Big( \sum_{j=0}^n \vartheta_{k,1,j,i',j'} \Big) c_{i',j'} = 0 \mbox{ for any } k\in \{1,\ldots,\om n-1\}
 \eneq
using the integer $\del_k$ (which depends also on $a$ and $n$) provided by Proposition~\ref{propgene}. This result asserts that~\eqref{eqsyslin} is a linear system of $M_0 = \om n-1$ equations in $N=a(n+1)$ unknowns $c_{i',j'} $, with integer coefficients bounded by 
\begin{equation} \label{eqet2}
\Big| \frac{\del_k}{(k-1)!} \sum_{j=0}^n \vartheta_{k,1,j,i',j'} \Big| \leq (n+1) k^a 2^n \del_k \leq \Big( 2(a+1)^{ \om}e^{3\om}\Big)^{ n (1+o(1))}
\eneq
as $n\to\infty$, since $k\leq \om n -1$ and $\om\geq 1$. To be consistent with the notation of Lemma~\ref{lemLS}, we let $H_k =\sqrt{N}(n+1)k^a 2^n \del_k $ for $1\leq k\leq M_0=\om n -1$.

\bigskip

In applying Lemma~\ref{lemLS}, 
for any $k\in\{\om n, \ldots, \gdom n-1\}$ we consider $\frakA_k$ given by Eq.~\eqref{eqdefad} as a linear combination of the unknowns $c_{i',j'} $, with integer coefficients bounded in absolute value by 
$ k^an^k$. We take $M = \gdom n-1$ and for each $k$ such that $M_0=\om n -1 <k \leq M$ we let $G_k = \sqrt N r^{\gdom n-k} $ and $H_k = \sqrt{N}k^an^k$. Then Lemma~\ref{lemLS} applies, and with its notation we have
$$X \leq \sqrt N \, \, \Big[ N^{(\gdom n-1)/2} \, \, \Big(2 (a+1)^{ \om}e^{3\om}\Big)^{(\om n -1)n (1+o(1))} \, \, \prod_{k=\om n}^{\gdom n -1} r^{\gdom n-k} \Big]^{\frac1{N-M_0}}$$
using Eq.~\eqref{eqet2}, so that
\begin{eqnarray*}
\log X 
&\leq& \frac{ n (1+o(1))}{a-\om} \Big( \om \log 2 + 3\om^2 + \om^2 \log (a+1) +\frac1{n^2} \sum_{k=\om n}^{\gdom n -1} (\gdom n-k) \log r\Big)\\
&\leq& \frac{ n (1+o(1))}{a-\om} \Big( \om \log 2 + 3\om^2 + \om^2 \log (a+1) +\frac12 \gdom^2 \log r\Big).
\end{eqnarray*}
This concludes the proof of Theorem~\ref{thsiegel}.

\section{Main part of the proof} \label{secdio}

In this section we prove Theorem~\ref{thintroun} stated in the introduction; we explain in \S \ref{subsecpolylogs} how to modify this proof and deduce Theorem~\ref{thpolylogs}.
 We explain the notation and sketch the proof in \S \ref{subsecdio1}. We obtain an expansion in polylogarithms in \S \ref{subsecdev}. Then we study the resulting linear forms: their coefficients (\S \ref{subsecdio3}) and their asymptotic behavior (\S \ref{subsecdio4}). We apply a multiplicity estimate in \S \ref{subseczero}, and conclude the proof in \S \ref{subsecfin}.

\subsection{Setting, notation and sketch of the proof} \label{subsecdio1}

Let $a, r , \om , \gdom \geq 1$ and $n\geq 2$, with $a,n\in\Z$, $r,\om, \gdom \in\Q$, and $1\leq \om\leq \gdom < a$; we assume $rn$, $\om n$ and $\gdom n$ to be integers. We shall use also another parameter $h\in \Z$, with $0\leq h \leq a$, to bound the order $p$ of derivation with respect to $t$. In our application, $a$, $r$, $\om$, $\gdom$, $h$ will be fixed and $n$ will tend to $\infty$. We refer to the end of this section (and to \S \ref{subsecfin}) for the choice of parameters.

 Using Siegel's lemma we have constructed in Theorem~\ref{thsiegel} (see \S \ref{subsec31}) integers $\cij\in\Z$, for $1\leq i \leq a$ and $0\leq j \leq n$, such that 
$$F_n(X) = \sum_{i=1}^a \sum_{j=0}^n \frac{\cij}{(X+j)^i} \in\Q(X)$$
satisfies $F_n(t) = O(|t|^{-\om n})$ as $\vert t \vert \to\infty$, with $|\cij| \leq \chi^{n(1+o(1))}$ as $n\to \infty$, where 
\begin{equation} \label{eqchi41}
\chi = \exp\Big(\frac{ \om \log 2 + 3\om^2 + \om^2 \log (a+1) +\frac12 \gdom^2 \log r }{a-\om }\Big).
\eneq
We have also 
\begin{equation} \label{eqmajofrakad}
|\frakA_d |\leq r^{d-\gdom n} n^d d^a \chi^{n(1+o(1))}
\eneq
for any $d < \gdom n$, where $ \frakA_d$ is defined by
\begin{equation} \label{eqfnad4}
 F_n(t) = \sum_{d=1}^\infty \frac{\frakA_d}{t^d} \, \, \mbox{ if $|t|>n$;}
\eneq
notice that the upper bound~\eqref{eqmajofrakad} is interesting only when $\om n \leq d < \gdom n$ since $\frakA_d = 0$ for any $d< \om n$.

\bigskip

For any $p\geq0$, the $p$-th derivative of $F_n$ is
$$F_n^{(p)}(X) = \sum_{i=1}^a \sum_{j=0}^n \frac{\cij (-1)^p (i)_p}{(X+j)^{i+p}} $$
where $(i)_p=i(i+1)\ldots (i+p-1)$. As mentioned at the beginning of this section, we 
 fix an additional parameter $h\geq 0$ with $h\leq a$.
 For any $z\in\C$ such that $\vert z \vert = 1$ and any $p\in\zeroh$ we consider 
$$
S_{n,p}(z) 
 = z^{rn} \sum_{t=rn+1}^\infty \Big( F_n^{(p)}(t) z^{-t} - F_n^{(p)}(-t) z^{ t} \Big)
$$
 which is convergent since $F_n^{(p)}(t) = O(|t|^{-\om n})$ as $\vert t \vert \to\infty$, with $\om n \geq 2$. 
The point here is that even zeta values should not appear in the linear combination we are trying to construct. A symmetry phenomenon (related to well-poised hypergeometric series) is used in general to obtain this property. However we have to consider derivatives of $S_{n,p}(z) $ to apply the multiplicity estimate, and this property is not transfered to derivatives. We overcome this difficulty as in~\cite{SFcaract}, by considering the functions $\Li_i(1/z)-(-1)^i \Li_i(z)$ instead of just $\Li_i(1/z)$. This leads to the definition above of $S_{n,p}(z) $, instead of simply $z^{rn} \sum_{t=rn+1}^\infty F_n^{(p)}(t) z^{-t}$. 
 
\bigskip
 
 We let also 
\begin{equation} \label{eq44bis}
P_i(z) = \sum_{j=0}^n \cij z^j \mbox{ for } 1\leq i \leq a
\eneq
and we shall prove in Lemma~\ref{lemdev} that, if $z\neq 1$, 
\begin{equation} \label{eq410}
S_{n,p}(z) 
 = V_p(z) + \sum_{i=1}^a z^{rn} P_i(z) (-1)^p (i)_p \Big( \Li_{i+p}(1/z) - (-1)^{i+p} \Li_{i+p}(z) \Big)
 \end{equation} 
for some polynomial $V_p\in\Q[X] $ of degree at most $2rn$. For $k\geq 1$ we shall consider the $(k-1)$-th derivative $S_{n,p}^{(k-1)}(z)$ of $S_{n,p}(z) $. Since the coefficients of the polynomial $V_p$ have large denominators (that would ruin our Diophantine application), we shall be interested only in integers $k$ such that $k-1\geq 2rn+1>\deg V_p$, so that $V_p^{(k-1)}=0$.

For $0\leq p \leq h$ and $1 \leq i \leq a$ we let 
 \begin{equation} \label{eqdefqp}
\qqp_{i+p}(z) = z^{rn} P_i(z) (-1)^p (i)_p
\eneq
and also $ \qqp_i(z)=0$ for $i \in \{1,\ldots,p\} \cup \{a+p+1,\ldots,\aplush\}$. Then Eq.~\eqref{eq410} reads
 \begin{equation} \label{eqdevsnp}
S_{n,p}(z)=V_p(z)+\sum_{i=1}^{\aplush} \qqp_i(z) \Big( \Li_i(1/z) - (-1)^{i} \Li_{i}(z)\Big).
\eneq
Now let $ \qqp_{1,0}(z) = 0$, $ \qqp_{1,i}(z) = \qqp_i(z)$ for any $i\in\{1,\ldots,\aplush\}$, and for $k\geq 2$:
 \begin{equation} \label{eqdefqpk}
\left\{ \begin{array}{l}
 \qqp_{k,i}(z) = { \qqppr_{k-1,i}}(z) - \frac1{z} \qqp_{k-1,i+1}(z) \mbox{ for } 1\leq i \leq \aplush \\
 \qqp_{k,0}(z) = { \qqppr_{k-1,0}}(z) + \frac{z+1}{z(1-z)} \qqp_{k-1,1}(z) 
\end{array}\right.
\eneq
where $ \qqp_{k-1,\aplush+1}$ is taken to be the zero function.  In particular we have $ \qqp_{k,i}(z)=0$ for any $i\in \{a+p+1,\ldots,\aplush\}$, but not (in general) for $0\leq i \leq p$. Since the derivative of $ \Li_i(1/z) - (-1)^{i} \Li_{i}(z) $ is $ \frac{z+1}{z(1-z)} $ for $i=1$, and $ - \frac1{z} \Big( \Li_{i-1}(1/z) - (-1)^{i-1} \Li_{i-1}(z)\Big)$ for $i\geq 2$, we have 
 \begin{equation} \label{eq33nv}
S_{n,p}^{(k-1)}(z)= \qqp_{k,0}(z) + \sum_{i=1}^{\aplush} \qqp_{k,i}(z) \Big( \Li_i(1/z) - (-1)^{i} \Li_{i}(z)\Big) \mbox{ for any } k\geq 2rn+2
\eneq
since $\deg V_p \leq 2rn$; when $1 \leq k \leq 2rn+1$ an additional term $V_p^{(k-1)}(z)$ appears on the right hand side. The point is that we have now many linear forms for each value of $n$, as $k$ and $p$ vary. This is necessary to apply the multiplicity estimate, and then Siegel's linear independence criterion. 

\bigskip

For any $k\geq 2rn+2$ we let 
 \begin{equation} \label{eqdeflikn}
\likn = (-2)^{k-1} \frac{\del_k}{(k-1)!} \qqp_{k,i} (-1) \mbox{ for } 0\leq i \leq \aplush 
\eneq
where $\del_k= \del_k (a+h,(r+1)n, 1,1)$ is given by Proposition~\ref{propgene} in \S \ref{subsectech1} with $a $ replaced by $\aplush$ and $n$ by $(r+1)n$; then Eq.~\eqref{eq33nv} yields
 \begin{equation} \label{eqFL}
 (-2)^{k-1} \frac{\del_k}{(k-1)!}
S_{n,p}^{(k-1)}(-1) = \lzkn + \sum_{i=1}^{\aplush} \likn ( 1- (-1)^{i}) \Li_i(-1) .
\eneq
These are the linear forms we are interested in, with $0\leq p \leq h$ and $2rn+2 \leq k\leq \capa n$ (where $\capa\in\Q$ is a fixed parameter such that $2r<\capa < \om$). We shall prove in Lemma~\ref{lemcoeffs} that their coefficients are not too large integers, namely 
$\likn \in\Z$ and
$$ | \likn | \leq \beta^{n(1+o(1))} \mbox{ with }
 \beta = \chi \Big( e^3 (2a+1) \Big)^\capa \cdot 4^{\capa+r+1}. $$
Then in Lemma~\ref{lemasy} we shall prove that these linear forms are small :
$$\Big| \lzkn + \sum_{i=1}^{\aplush} \likn \Big( 1 - (-1)^i\Big) \Li_i(-1) \Big| \leq \alpha^{n(1+o(1))} \mbox{ with } 
\alpha = \chi r^{-\gdom}( 2e^4(2a+1))^{\capa} . $$
Assume that $(h+1)(\capa-2r)+\om >a$, and that $n$ is sufficiently large.
Then using the generalization of Shidlovsky's lemma stated in \S \ref{subsecshidenonce} we prove
 in \S \ref{subseczero} that there are sufficiently many linearly independent linear forms among them; this allows us in \S \ref{subsecfin}
to apply Siegel's linear independence criterion 
(recalled in \S \ref{subsecsiegel}) and deduce that 
$$\dim_\Q \Span_\Q (\{1\}\cup \{ (1-(-1)^i) \Li_i(-1), \, 1\leq i \leq \aplush\}) \geq 1 - \frac{\log \alpha}{\log \beta } .$$
Choosing appropriate parameters (namely $r =3.9$, $\capa = 10.58$, $\om = 11.58$, $\gdom \in\Q$ sufficiently close to $3.9 \sqrt{a \log a}$, and $h =0.36\ a$) enables one to conclude the proof of Theorem~\ref{thintroun} (see \S \ref{subsecfin} for details); recall that $(1-(-1)^i)\Li_i(-1)$ vanishes when $i$ is even, and is equal to $2(2^{1-i}-1)\zeta(i)$ when $i\geq 3 $ is odd.

\subsection{Expansion in polylogarithms} \label{subsecdev}

\begin{Lem} \label{lemdev} 
For any $p\in\zeroh$ there exists a polynomial $V_p\in\Q[X] $ of degree at most $2rn$ such that, for any $z\in\C$ with $\vert z\vert =1$ and $z \neq 1$, 
$$
S_{n,p}(z) 
 = V_p(z) + \sum_{i=1}^a z^{rn} P_i(z) (-1)^p (i)_p \Big( \Li_{i+p}(1/z) - (-1)^{i+p} \Li_{i+p}(z) \Big).$$
\end{Lem}

\Dem of Lemma~\ref{lemdev}: 
 To begin with, we let 
 \begin{equation} \label{eqdefsninf} 
 \Snpinf(z) 
 = z^{rn} \sum_{t=rn+1}^\infty F_n^{(p)}(t) z^{-t}
 \eneq
 for $z\in \C$, $|z|\geq1$, $z\neq 1$. We have
\begin{eqnarray*} 
 \Snpinf(z) 
&=& \sum_{t=rn+1}^\infty \sum_{i=1}^a \sum_{j=0}^n \frac{\cij (-1)^p (i)_p}{(t+j)^{i+p}} z^{rn-t} \\
&=& \sum_{i=1}^a \sum_{j=0}^n \cij (-1)^p (i)_p \sum_{\ell =rn+1+j}^\infty \frac{z^{rn-\ell+j}}{ \ell^{i+p}} \\
&& \quad \quad \quad \mbox{since this series is convergent (because $|z|\geq1$ and $z\neq 1$)} \\
&=& \sum_{i=1}^a \sum_{j=0}^n \cij (-1)^p (i)_p \Big( z^{rn+j} \Li_{i+p}(1/z) - \sum_{\ell =1}^{rn+j} \frac{z^{rn-\ell+j}}{ \ell^{i+p}}\Big) 
\end{eqnarray*}
so that
$$
\Snpinf(z) = \Vpinf(z) + \sum_{i=1}^a z^{rn} P_i(z) (-1)^p (i)_p \Li_{i+p}(1/z)$$
where (as defined above) 
$$P_i(z) = \sum_{j=0}^n \cij z^j \mbox{ for } 1\leq i \leq a$$
and 
\begin{equation} \label{eqdefVzero}
\Vpinf(z) = - \sum_{i=1}^a \sum_{j=0}^n \cij (-1)^p (i)_p \sum_{t=0}^{rn+j-1} \frac{z^{t}}{(rn+j-t)^{i+p}} \in\Q[z].
\eneq
Observe that the polynomials $P_i$ have degree at most $n$, and do not depend on $p$, whereas $\Vpinf$ depends on $p$ and has degree at most $(r+1)n-1$.

\bigskip

On the other hand we consider, for $z\in \C$ with $|z|\leq1$ and $z\neq 1$, 
\begin{eqnarray*} 
\Snpzero(z) 
&=& z^{rn} \sum_{t=rn+1}^\infty F_n^{(p)}(-t) z^{t} 
\\
&=& \sum_{t=rn+1}^\infty \sum_{i=1}^a \sum_{j=0}^n \frac{\cij (-1)^p (i)_p}{(-t+j)^{i+p}} z^{rn+t} \\
&=& \sum_{i=1}^a \sum_{j=0}^n \cij (-1)^p (i)_p (-1)^{i+p} \sum_{\ell =rn+1-j}^\infty \frac{z^{rn+\ell+j}}{ \ell^{i+p}}\nonumber \\
&=& \sum_{i=1}^a \sum_{j=0}^n \cij (-1)^p (i)_p (-1)^{i+p} \Big( z^{rn+j} \Li_{i+p}(z) - \sum_{\ell =1}^{rn-j} \frac{z^{rn+\ell+j}}{ \ell^{i+p}}\Big)\nonumber 
\end{eqnarray*}
so that
$$
\Snpzero(z) = \Vpzero(z) + \sum_{i=1}^a z^{rn} P_i(z) (-1)^p (i)_p (-1)^{i+p} \Li_{i+p}(z)$$
with the same polynomials $P_i$, and 
\begin{equation} \label{eqdefVinf}
\Vpzero(z) = - \sum_{i=1}^a \sum_{j=0}^n \cij (-1)^i (i)_p \sum_{t=rn+j+1}^{2rn} \frac{z^{t}}{(t-rn -j)^{i+p}} \in\Q[z].
\eneq
Observe that $\Vpzero$ has degree at most $2rn$ and is a multiple of $z^{rn+1}$. Since $S_{n,p}(z) = \Snpinf(z) - \Snpzero(z) $, we let $V_p(z) = \Vpinf(z) - \Vpzero(z) $; this concludes the proof of Lemma~\ref{lemdev}.

\subsection{Coefficients of the linear forms} \label{subsecdio3}

For any algebraic number $\xi$, we denote by $\house{\xi}$ its house, i.e. the maximum modulus of its Galois conjugates. To prepare for the proof of Theorem~\ref{thpolylogs} (see \S \ref{subsecpolylogs}) we shall estimate the coefficients of the linear forms in a slightly more general setting than what is needed in the proof of Theorem~\ref{thintroun}.

Let $z_0\in\Qbar$ be such that $|z_0| \geq 1$ and $z_0\neq 1$; denote by $\qq\in\N\etoile$ a denominator of $z_0$, i.e. such that $\qq z_0 \in \OQzz$ where $\OQzz$ is the ring of integers of $\Q(z_0)$.
For any $k\geq 1 $ we let 
 \begin{equation} \label{eqdefl43}
\likn(z_0) = q^{(r+1)n+k-1} z_0^{k-1}(1-z_0)^{k-1} \frac{\del_k}{(k-1)!} \qqp_{k,i} (z_0) \mbox{ for } 0\leq i \leq \aplush 
\eneq
where $\del_k=\del_k (a+h,(r+1)n,1,1)$ is given by Proposition~\ref{propgene} in \S \ref{subsectech1}, and the rational functions $ \qqp_{k,i} (z)$ are defined by Eq.~\eqref{eqdefqpk}. 
The special case needed in the proof of Theorem~\ref{thintroun} is $z_0=-1$, $\qq=1$; then $\Q(z_0) = \Q$ and $\OQzz = \Z$, and 
$ \likn(z_0) = \likn$ (see Eq.~\eqref{eqdeflikn}).

\begin{Lem} \label{lemcoeffs} 
We have $\likn(z_0) \in\OQzz$ for any $p\in\zeroh$, any $i\in\{0,\ldots,\aplush\}$ and any $k\geq 1$. Moreover, provided $k\leq\capa n$ with a fixed $\capa\geq r+1$ (independent from $n$), we have as $n\to\infty$:
$$ \house{\likn(z_0)} \leq \beta^{n(1+o(1))} \mbox{ with } \beta = \chi \Big( 8 e^3 (2a+1) \Big)^\capa \cdot \Big( \qq \max(1, \house{z_0} , \house{ 1-z_0} ) \Big)^{\capa+r+1} $$
where $\chi$ is defined by Eq.~\eqref{eqchi41}.
\end{Lem}

\Dem of Lemma~\ref{lemcoeffs}: 
We fix $p$ and apply the results of \S \ref{subsectech1}. With respect to the notation of that section, $P_i(z)$ is replaced with $\qqp_i(z)$, $a$ with $\aplush$ and $n$ with $(r+1)n$; recall that $\deg \qqp_i\leq (r+1)n$ for any $i\in\{1,\ldots,\aplush\}$ (see Eq.~\eqref{eqdefqp} and the line following it). We take $\bte=\ate = 1 $ in the notation of \S \ref{subsectech1}, so that Eqns.~\eqref{eqrecpkijtechnique} and~\eqref{eqdefqpk} are consistent.
We write 
$$
\left\{ \begin{array}{l}
 z^{k-1} \qqp_{k,i}(z) = \sum_{j=0}^{(r+1)n} q_{k,i,j} z^j \mbox{ if } i\geq 1,\\
z^{k-1} (1-z)^{k-1} \qqp_{k,0}(z)= \sum_{j=0}^{(r+1)n+k-1} q_{k,0,j} z^j.
\end{array}\right.
$$
Then Eq.~\eqref{eqdefl43} reads 
 \begin{equation} \label{eq792}
 \likn(z_0) = 
\qq^{k-1} (1-z_0)^{k-1} \sum_{j=0}^{(r+1)n} \frac{\del_k}{(k-1)!} q_{k,i,j} \qq^{(r+1)n} z_0^j \mbox{ for } 1\leq i \leq \aplush,
\eneq
and
 \begin{equation} \label{eq793}
\lzkn(z_0) = \sum_{j=0}^{(r+1)n+k-1} \frac{\del_k}{(k-1)!} q_{k,0,j} \qq^{(r+1)n+k-1} z_0^j.
\eneq 
To be consistent with the notation of \S \ref{subsectech1} we write also $\qqp_i(z)= \sum_{j=0}^{(r+1)n} c'_{i,j} z^j$ for $1\leq i \leq \aplush$. Combining Eq.~\eqref{eqch0} with part $(ii)$ of Proposition~\ref{propgene}, we deduce that $ \frac{\del_k}{(k-1)!} q_{k,i,j} \in\Z$ for any $k$, $i$, $j$, since $ c'_{i',j'} \in\Z$ for any $i'$, $j'$. Moreover, part $(iii)$ of Proposition~\ref{propgene} and Eq.~\eqref{eqch0} yield
$$ \Big| \frac{\del_k}{(k-1)!} q_{k,i,j} \Big| \leq k^{a+h+1} \, \, 8^{\max(k, (r+1)n)} \, \,\del_k \, \, (a+h) \, \,((r+1)n+1) \, \,\max_{i',j'} |c'_{i',j'}| $$
for any $k$, $i$, $j$, with $ \del_k \leq ( e^3 (a+h+1) )^{\max(k, (r+1)n)}$ according to part $(i)$ -- recall that Proposition~\ref{propgene} is applied with $\aplush $ and $ (r+1)n$ instead of $a$ and $n$, respectively. Since $a+h\leq 2a$, we deduce that 
$$ \Big| \frac{\del_k}{(k-1)!} q_{k,i,j} \Big| \leq 2 \, \,k^{2a+1} \, \, (8 e^3 (2a+1) )^{\max(k, (r+1)n)} \, \,a((r+1)n+1) \, \,\max_{i',j'} |c'_{i',j'}|. $$
Using Eqns.~\eqref{eq792} and~\eqref{eq793} we obtain $ \likn(z_0) \in \OQzz$ for any $i \in \{0,\ldots,2a\}$, any $k\geq 1$ and any $p\in\zeroh$, and 
\begin{eqnarray*}
 \vert \likn(z_0) \vert &\leq &
2 \, \, k^{2a+1} \, \,( 8 e^3 (2a+1) )^{\max(k, (r+1)n)} \, \,a((r+1)n+k)^2 \, \,\max_{i',j'} |c'_{i',j'}| \\
 &&\quad \cdot \, 
\qq^{(r+1)n+k-1} \, \, \max(1, \house{ z_0}^{(r+1)n}) \, \,\max( 1, \house{ 1-z_0} ^{k-1} , \house{z_0}^{k-1} ) .
\end{eqnarray*}
Now Eq.~\eqref{eqdefqp} and Theorem~\ref{thsiegel} yield $ \max_{i',j'} |c'_{i',j'}| \leq (a)_a \chi^{n(1+o(1))}$ since $h\leq a$. Using the assumption $k\leq \capa n$ with $\capa\geq r+1$, this concludes the proof of Lemma~\ref{lemcoeffs}.

\subsection{Asymptotic estimate of the linear forms} \label{subsecdio4}

 Let $z_0\in\Qbar$ be such that $|z_0| = 1$; in this section $z_0$ could be equal to 1. We shall take $z_0=-1$ in the proof of Theorem~\ref{thintroun}, and adapt the proof of Lemma~\ref{lemasy} below in \S \ref{subsecpolylogs} to prove Theorem~\ref{thpolylogs}. Recall that $\del_k=\del_k (a+h,(r+1)n, \bte,\ate)\in\N\etoile$ has been defined in Proposition~\ref{propgene} (in which $a $ should be replaced with $\aplush$ and $n$ by $(r+1)n$), and $\chi$ in Theorem~\ref{thsiegel}.

\begin{Lem} \label{lemasy}
Assume that $r\geq 2$, $0 \leq p \leq h$, and $2rn+2\leq k \leq \capa n $, with $ \capa < \om$. Then we have
$$\Big| \frac{\del_k}{(k-1)!} S_{n,p}^{(k-1)}(z_0)\Big| \leq \alpha_0^{n(1+o(1))} \mbox{ with } 
\alpha_0 = \chi r^{-\gdom}( e^4(2a+1))^{\capa} .$$
\end{Lem}

\Dem of Lemma~\ref{lemasy}: 
Recall that $S_{n,p}(z)=\Snpinf(z)- \Snpzero(z)$ with the notation introduced in the proof of Lemma~\ref{lemdev}. 
Taking the $p$-th derivative of Eq.~\eqref{eqfnad4} (see \S \ref{subsecdio1}) yields 
$F_n^{(p)}(t)=\sum_{d=1}^{\infty} \frac{\frakA_d (-1)^p (d)_p}{t^{d+p}} $
 for $|t|>n$. By definition of $\Snpinf(z) $ (see Eq.~\eqref{eqdefsninf} in \S \ref{subsecdev}) we obtain
 \begin{equation} \label{eq1235}
\Snpinf(z) = \sum_{t=rn+1}^\infty \sum_{d=1}^\infty \frac{\frakA_d (-1)^p (d)_p}{t^{d+p}} z^{rn-t} \mbox{ for } |z| \geq 1, \, z\neq 1.
\eneq
Now Theorem~\ref{thsiegel} asserts that $F_n(t)=O(|t|^{-\om n})$ as $|t|\to\infty$, so that 
 $\frakA_d=0$ for any $d\in\{1,\ldots,\om n -1\}$: the sum on $d$ in Eq.~\eqref{eq1235} starts only at $d=\om n$. 
 Therefore we have for any $k\geq1$:
$$ \frac{\del_k}{(k-1)!} \Snpinfkmu (z) = (-1)^{k-1} \del_k \sum_{t=rn+1}^\infty \sum_{d= \om n}^\infty \frac{\frakA_d (-1)^p (d)_p}{t^{d+p}} \binom{t-rn+k-2}{k-1} z^{rn-t-k+1}.$$
Since $\vert z \vert \geq 1$ and $t^p \geq 1$ we obtain 
$$ \Big| \frac{\del_k}{(k-1)!} \Snpinfkmu(z)\Big| \leq \del_k \sum_{t=rn+1}^\infty \binom{t-rn+k-2}{k-1} 
 \Big( \frac{n}{t }\Big)^{\om n}
\sum_{d= \om n}^\infty \frac{ | \frakA_d | (d)_p }{t^{d-\om n}}n^{-\om n}. $$
We bound $ | \frakA_d | $ trivially (using Eq.~\eqref{eqdefad} and assertion $(ii)$ of Theorem~\ref{thsiegel}) for $d\geq \gdom n$, and we use assertion $(iii)$ of Theorem~\ref{thsiegel} for $d$ such that $\om n \leq d < \gdom n$. Therefore we have
 \begin{equation} \label{eq12357}
\Big| \frac{\del_k}{(k-1)!} \Snpinfkmu(z)\Big| \leq \del_k \chi^{n(1+o(1))} \sum_{t=rn+1}^\infty \binom{t-rn+k-2}{k-1} 
 \Big( \frac{n}{t }\Big)^{\om n}
\sum_{d= \om n}^\infty u_{t,d}
\eneq
where the sequence $o(1)$ does not depend on $k$, nor on $p$, and tends to 0 as $n\to\infty$; we define $u_{t,d}$ by 
$$u_{t,d} = (d)_p d^a (n/t)^{d-\om n} \mbox{ for } d \geq \gdom n$$
and 
$$u_{t,d} = r^{d-\gdom n} (d)_p d^a (n/t)^{d-\om n} \mbox{ for } \om n \leq d < \gdom n.$$
Let us bound the term $\sum_{d= \om n}^\infty u_{t,d}$ in Eq.~\eqref{eq12357}. For any $d \geq \gdom n$ we have $u_{t,d+1}/u_{t, d} \leq (1+\frac{p}{d}) \cdot (1+\frac1{d})^a \cdot \frac1r \leq\frac3{2r}$ for any $t\geq rn+1$, provided $n$ is large enough (using the assumption that $\gdom>0$). Since $r\geq 2$ we obtain
\begin{equation} \label{eq319ter}
\sum_{d= \gdom n}^\infty u_{t,d} \leq u_{t,\gdom n} \sum_{d= \gdom n}^\infty \Big(\frac34\Big) ^{d-\gdom n} \leq 4 r^{(\om-\gdom) n}(\gdom n)_p (\gdom n)^a 
\end{equation}
for any $t\geq rn+1$. On the other hand, for $\om n \leq d < \gdom n$ we have 
$$u_{t,d} = r^{(\om-\gdom) n} (d)_p d^a (rn/t)^{d-\om n} \leq r^{(\om-\gdom) n} ( \gdom n)_p ( \gdom n)^a .$$
Combining this upper bound with Eq.~\eqref{eq319ter} yields
$$
\sum_{d= \om n}^\infty u_{t,d} \leq (4 + (\gdom-\om)n) r^{(\om-\gdom) n} (\gdom n)_p (\gdom n)^a \leq r^{(\om-\gdom) n} \chi^{o(n)}; 
$$
here and below, the sequences $o(\cdots)$ may depend on $p$ (but not on $k$). Using Eq.~\eqref{eq12357} we obtain
\begin{equation} \label{eq319bis}
 \Big| \frac{\del_k}{(k-1)!} \Snpinfkmu (z)\Big| \leq r^{-\gdom n} \del_k \chi^{n(1+o(1))} \sum_{t=rn+1}^\infty \binom{t-rn+k-2}{k-1} \Big( \frac{rn}{t }\Big)^{\om n} .
\eneq
We let $\sigma = \frac{k-1}{rn}$ so that $\sigma > 1$. Let $t>rn$; then we have $ t-rn+k-2 \leq t + (\sigma-1)rn < \sigma t$ so that 
$$ \binom{t-rn+k-2}{k-1} \Big( \frac{rn}{t }\Big)^{\om n-2} \leq \frac{(\sigma t)^{k-1}}{(k-1)!} \Big( \frac{rn}{t }\Big)^{\om n-2} \leq
 \frac{\sigma ^{k-1} (rn) ^{k-1} }{(k-1) ^{k-1} e^{-k+1}} \Big( \frac{rn}{t }\Big)^{\om n-k-1} \leq e^{k-1}$$
 since $\frac{rn}{t }\leq 1$ and $k+1\leq \capa n +1 \leq \om n$; recall that $(k-1)!\geq (\frac{k-1}{e})^{k-1}$, 
and $\sigma r n = k-1$ by definition of $\sigma$. This proves that 
 \begin{equation} \label{eqmajolemasytemp}
 \sum_{t=rn+1}^\infty \binom{t-rn+k-2}{k-1} \Big( \frac{rn}{t }\Big)^{\om n } \leq r^2 n^2 e^{k-1} \pi^2/6.
\eneq
Using Eq.~\eqref{eq319bis}, Theorem~\ref{thsiegel} and assertion $(i)$ of Proposition~\ref{propgene} (where $a$ is replaced with $\aplush\leq 2a$ and $n$ with $(r+1)n$), we obtain 
$$ \Big| \frac{\del_k}{(k-1)!} \Snpinfkmu (z)\Big| \leq \alpha_0^{n(1+o(1))} .$$

\bigskip

We now turn to $ \Snpzerokmu (z)$ (recall that $S_{n,p}(z) = \Snpinf (z) - \Snpzero (z)$). As for $\Snpinf$ above, we have 
$$\Snpzero(z) = \sum_{t=rn+1}^\infty \sum_{d=\om n }^\infty \frac{\frakA_d (-1)^p (d)_p}{(-t)^{d+p}} z^{rn+t} \mbox{ for } |z| \leq 1, \, z\neq 1,$$
so that, for any $k\geq 2rn+2$, 
$$ \frac{\del_k}{(k-1)!} \Snpzerokmu (z) = \del_k \sum_{t=k-1-rn}^\infty \sum_{d= \om n}^\infty \frac{\frakA_d (-1)^d (d)_p}{t^{d+p}} \binom{rn+t}{k-1} z^{rn+t-k+1}.$$
We have 
$$ \Big| \frac{\del_k}{(k-1)!} \Snpzerokmu (z) \Big| \leq \del_k \chi^{n(1+o(1))} \sum_{t=k-1-rn}^\infty \binom{rn+t}{k-1} 
 \Big( \frac{n}{t }\Big)^{\om n}
\sum_{d= \om n}^\infty u_{t,d}$$
with the same $ u_{t,d}$ as above, so that 
\begin{equation} \label{eq319bisdeux}
 \Big| \frac{\del_k}{(k-1)!} \Snpzerokmu (z) \Big| \leq r^{-\gdom n } \del_k \chi^{n(1+o(1))} \sum_{t=k-1-rn}^\infty \binom{rn+t}{k-1} \Big( \frac{rn}{t }\Big)^{\om n} . 
\eneq
Now we have $t+rn\leq t-rn+k-2$ for any $t$, so that $ \binom{rn+t}{k-1} \leq \binom{t-rn+k-2}{k-1} $: we obtain the same upper bound as in Eq.~\eqref{eq319bis}, and deduce in the same way
 $$ \Big| \frac{\del_k}{(k-1)!} \Snpzerokmu (z) \Big| \leq \alpha_0^{n(1+o(1))} .$$
 Since $S_{n,p}^{(k-1)}(z ) = \Snpinfkmu (z) - \Snpzerokmu (z) $, taking $z=z_0$ this concludes the proof of Lemma~\ref{lemasy}.

\subsection{Multiplicity estimate} \label{subseczero}

In this section we apply the multiplicity estimate stated in \S \ref{subsecshidenonce} to prove Proposition~\ref{proplemzero} below, which makes it possible to  apply the refinement of  Siegel's linear independence criterion proved in \S \ref{subsecsiegel}. 

\bigskip

To state Proposition~\ref{proplemzero}, recall that $P_i(z) = \sum_{j=0}^n \cij z^j$ for $1\leq i \leq a$. Since the integers $\cij$ are not all zero, we may consider 
$$b = \max\{ i\in\{1,\ldots,a\}, \, \exists j \in\zeron, \, \cij\neq 0\}.$$
Then we have $1\leq b \leq a$, $P_b\neq 0$, and $P_{b+1}= \ldots = P_a=0$. 
Eqns.~\eqref{eqdefqp},~\eqref{eqdefqpk} and~\eqref{eqdeflikn} show that $\qqp_i(z)$, $\qqp_{k,i}(z)$ and $\likn$ all vanish when $b+p+1\leq i \leq \aplush$: Eq.~\eqref{eqFL} becomes a linear form in 1 and the numbers $(1-(-1)^i) \Li_{i}(-1) $ for $1 \leq i \leq \bplush$, namely
 \begin{equation} \label{eqformeslinlemz}
 (-2)^{k-1} \frac{\del_k}{(k-1)!}
S_{n,p}^{(k-1)}(-1) = \lzkn + \sum_{i=1}^{\bplush} \likn ( 1- (-1)^{i}) \Li_i(-1) 
\eneq
 with $2rn+2 \leq k\leq \capa n$ and $0\leq p \leq h$. To sum up, for a given $n$ we have (small) linear forms,  indexed by $k$ and $p$, in  $\bplush+1$ numbers.
 
\bigskip

Usually, the conclusion of a zero estimate in this setting would be that there exist  $\bplush+1$ linearly independent linear forms among them. However this is {\em false} in general in our setting (see Remark \ref{remmatnoninv} below): there may be non-trivial linear relations between the coefficients $ \likn $, $0\leq i \leq \bplush$, valid for any $k$ and any $p$. The crucial point is that such a relation   cannot involve $\lzkn$, as the following result shows. This is sufficient to apply the refinement of  Siegel's linear independence criterion proved in \S \ref{subsecsiegel}.

\begin{Prop} \label{proplemzero} 
Assume that $(h+1)(\capa-2r)+\om > a$, and that $n$ is sufficiently large.  Let $x_0, \ldots, x_{\bplush}\in\Qbar$ be such that 
$$
\sum_{i=0}^{\bplush}  \likn  x_i =0 \mbox{ for any } k \in \intk \mbox{ and any } p \in \zeroh.
$$
Then $x_0=0$.
\end{Prop}

\begin{Remark} \label{remmatnoninv} 
Since our linear forms are constructed using Siegel's lemma, it seems extremely difficult to exclude the case where, 
  for some $\lambda\in\C$, we would have 
\begin{equation}\label{eqdrame}
\sum_{i=1}^b  Q^{[0]}_i(z) \frac{  (\lambda-\log z)^{i-1}}{(i-1)!} = O((z+1)^{\capa n})
\eneq
as $z\to -1$. Indeed, if $\lambda-\log (-1)\in\Z$, where we fix a determination of $\log z$ around $z=-1$, this amounts to $\capa n$ linear equations in the coefficients $\cij$; recall that these $a(n+1)$ coefficients have been constructed in \S \ref{subsecdio2} by solving $\om n -1$ linear equations, and $\om n -1 + \capa n$ is still much smaller that $a(n+1)$ with the parameters we shall choose in \S \ref{subsecfin}.

In case Eq. \eqref{eqdrame} holds, we deduce using Eq. \eqref{eqdefqp} that  for any $p\in\zeroh$,  
\begin{eqnarray}
\sum_{i=p+1}^{p+b}  Q^{[p]}_i(z) \frac{  (\lambda-\log z)^{i-1}}{(i-1)!}  
&=& \sum_{i=p+1}^{p+b}  Q^{[0]}_{i-p}(z) (-1)^p (i-p)_p \frac{  (\lambda-\log z)^{i-1}}{(i-p)_p(i-p-1)!}    \label{eqfctdrame}  \\
&=& (-1)^p  (\lambda-\log z)^p \sum_{j=1}^b  Q^{[0]}_j(z) \frac{  (\lambda-\log z)^{j-1}}{(j-1)!}  \nonumber  \\
&=& O((z+1)^{\capa n})  \mbox{ as $z \to -1$.}  \nonumber 
\end{eqnarray}
Recall that the $\qqp_{k,i}(z)$ were defined in Eq.~\eqref{eqdefqpk}  to compute derivatives of linear forms in the functions 
1 and $\Li_{i}(1/z )- (-1)^i \Li_{i}( z )$, $1 \leq i \leq \bplush$. Now the functions $ \frac{  (\lambda-\log z)^{i-1}}{(i-1)!}$ 
satisfy the same rules of differentiation so that $\sum_{i=p+1}^{p+b}  Q^{[p]}_{k,i}(z)  \frac{  (\lambda-\log z)^{i-1}}{(i-1)!}  $ is the $(k-1)$-th derivative of the function \eqref{eqfctdrame}: it vanishes at $z=-1$ for any  $k \in \intk $ and any $ p \in \zeroh$.  Using Eq.~\eqref{eqdeflikn} we obtain a non-trivial linear relation, valid for any $k$ and any $p$, between the coefficients $\likn $ of our linear forms: under the assumptions of Proposition \ref{proplemzero}, it would be false  to claim that $x_0=\ldots=x_{\bplush}=0$.
\end{Remark}

\begin{Remark} 
Let us comment on the assumption $(h+1)(\capa-2r)+\om > a$. To explain how necessary it is, we claim that if $(h+1)(\capa-2r)+\om < a$ then our approach cannot even exclude the case where $ (1-(-1)^i) \Li_{i}(-1) \in\Q$ for any $1 \leq i \leq \aplush$. The point is that the coefficients $\cij$ are provided by Siegel's lemma: they are not explicit, and the only property we can reasonably use in a multiplicity estimate is that $F_n(t)=O(t^{-\om n})$ as $|t|\to \infty$ (see Theorem~\ref{thsiegel}). This amounts to $\om n +O(1)$ linear equations in the unknowns $\cij$, where $O(1)$ denotes a term that is bounded uniformly with respect to $n$. Assuming that $ (1-(-1)^i) \Li_{i}(-1) \in\Q$ for any $1 \leq i \leq \aplush$, we claim that all linear forms~\eqref{eqformeslinlemz} may vanish, for any $2rn+2 \leq k\leq \capa n$ and any $0\leq p \leq h$. Indeed this would mean that the integers $\cij$ are solution of a linear system of $(h+1)(\capa-2r)n+\om n +O(1)$ linear equations with rational coefficients
(see Eqns.~\eqref{eqdeflikn},~\eqref{eqdefqp} and~\eqref{eq44bis}). If $(h+1)(\capa-2r)+\om < a$ and $n$ is sufficiently large, this system has fewer equations that the number of unknowns $\cij$ (namely, $a(n+1)$): there is a family of integers $\cij$, not all zero, that satisfy these equations. We see no reasonable way to prove that Theorem~\ref{thsiegel} does not provide this family; and if it does, all linear forms we are interested in vanish. Therefore we cannot hope to reach any contradiction if $(h+1)(\capa-2r)+\om < a$. 
\end{Remark}

In this section we prove Proposition~\ref{proplemzero}. To get ready for \S \ref{subsecpolylogs} (where the proof of Theorem~\ref{thintroun} is adapted to prove Theorem~\ref{thpolylogs}), we consider any $z_0\in\Qbar \setminus \{0,1\}$, not only the special case 
 $z_0=-1$ used to prove Theorem~\ref{thintroun}. In this general setting, the coefficients $  \likn  $ are defined by Eq.~\eqref{eqdefl43}.

\bigskip

Let $x_0, \ldots, x_{\bplush}\in\Qbar$ be as in Proposition~\ref{proplemzero}. By contradiction we assume $x_0\neq 0$, and even $x_0=1$ (dividing all $x_i$ by $x_0$ if necessary).  Using Eq.~\eqref{eqdefl43} we obtain
\begin{equation}\label{eqannulQx}
\sum_{i=0}^{\bplush} \qqp_{k,i}(z_0) x_i =0 \mbox{ for any } k \in \intk \mbox{ and any } p \in \zeroh.
\eneq
Throughout the proof of Proposition~\ref{proplemzero} we fix a small open disk centered at $z_0$, contained in $\C\setminus\{0,1\}$; all functions of $z$ we consider will be holomorphic on this disk.

We denote by $Y'=A_0 Y$ with $A_0\in M_{\bplush+1}(\Q(z))$ the differential system satisfied by the vector $Y(z)=\tra(y_0(z),\ldots,y_{b+h}(z))$ given by $y_0(z)=1$ and $y_i(z)= \Li_{i}(1/z )- (-1)^i \Li_{i}( z )$ for $1 \leq i \leq \bplush$. Since $z_0\not\in \{0,1\}$, the point $z_0$ is not a singularity of this system: there exists a solution $\tra(g_0(z),\ldots,g_{b+h}(z))$ of this system consisting in functions holomorphic around $z_0$ such that $g_i(z_0)=x_i $ for any $0\leq i \leq \bplush$. We have
$$
 g'_0(z)=0, \quad g'_1(z)=\frac{z+1}{z(1-z)} g_0(z), \quad \mbox{ and  } g'_i(z)=\frac{-1}{z} g_{i-1}(z) \mbox{ for  } 2\leq i \leq \bplush.$$
We consider, for any $p\in \zeroh$, the function
\begin{equation}\label{eqdefff}
\ff_p(z) = \TT_p(z) + \sum_{i=0}^{\bplush} \qqp_{i}(z ) g_i(z)
\eneq
where $\TT_p(z)\in\Qbar[z]_{\leq 2rn }$ is chosen so that $\ff_p(z) = O((z-z_0)^{2rn+1})$ as $z\to z_0$ (namely, $- \TT_p(z)$ is the Taylor approximation polynomial of degree at most $2rn $ of $ \sum_{i=0}^{\bplush} \qqp_{i}(z ) g_i(z)$ around $z_0$).

\bigskip

\noindent {\bf Step 1:} Vanishing of $\ff_p(z) $ with order at least $\capa n$ at $z_0$.

We claim that for any $p\in\zeroh$ we have
\begin{equation}\label{eqannulff}
\ff_p(z) = O((z-z_0)^{ \capa n})\mbox{ as } z\to z_0.
\eneq
Indeed the definition of $\qqp_{k,i}(z)$ in Eq.~\eqref{eqdefqpk}, intended to compute derivatives of linear forms in the functions 
1 and $\Li_{i}(1/z )- (-1)^i \Li_{i}( z )$, $1 \leq i \leq \bplush$ (see Eq.~\eqref{eqdevsnp}), can also be used for linear forms in $g_0(z)$, \ldots, $g_{\bplush}(z)$ because they satisfy the same rules of differentiation (i.e., $\tra(g_0(z),\ldots,g_{b+h}(z))$ is a solution of $Y'=A_0 Y$). Therefore we have
$$
\ff_p ^{(k-1)}(z) = \TT_p ^{(k-1)}(z)+ \sum_{i=0}^{\bplush} \qqp_{k,i}(z ) g_i(z) \mbox{ for any } k \geq 1.
$$
For any $k\in \intk$, Eq.~\eqref{eqannulQx} yields $\ff_p ^{(k-1)}(z_0)=0$ since $g_i(z_0)=x_i$ and $\deg \TT_p\leq 2rn $. This concludes the proof of Eq.~\eqref{eqannulff}.

\bigskip

\noindent {\bf Step 2:} Defining new polynomials and functions.

The strategy of the proof of Proposition~\ref{proplemzero} is to apply Shidlovsky's lemma. The problem for now is that the functions $\ff_p$ are not suitable for this: the polynomials $\qqp_{i}(z ) $ in Eq.~\eqref{eqdefff} should be independent from $p$. Their dependence in $p$ is rather weak (see Eq.~\eqref{eqdefqp}), and we shall overcome this difficulty now (see Eqns.~\eqref{eqcclhh} and~\eqref{eqsysdiff}). 

We consider the functions $\hh_q(z)$ defined by:
\begin{equation}\label{eqdefhh}
\hh_q(z) = \sum_{p=0}^q \bino (-\log z)^{q-p} \ff_p(z) \mbox{ for } q \in\zeroh;
\eneq
here and throughout \S \ref{subseczero}, $\log z$ can be seen formally.
We define also $\yy_{0,q},\ldots, \yy_{\bplush,q}$ for $q\in\zeroh$ by:
\begin{equation}\label{eqdefyq}
\left\{
\begin{array}{rcl}
\yy_{i,q}(z) &=&
0 \mbox{ for } 0\leq i \leq h-q-1 \\ \\
\yy_{i,q}(z) &=&
\frac{q!}{(i+q-h)!} (-\log z)^{i+q-h} \mbox{ for } h-q \leq i \leq h \\ \\
\yy_{i,q}(z) &=&
 \sum_{p=0}^q \bino (-\log z)^{q-p} (-1)^p (i-h)_p g_{i-h+p}(z) \mbox{ for } h+1 \leq i \leq \bplush
 \end{array}
 \right.
\eneq
and the following polynomials $\SS_0,\ldots,\SS_{\bplush}\in\Qbar[z]_{\leq 2rn}$:
\begin{equation}\label{eqdefSS}
\left\{
\begin{array}{rcl}
S_i(z) &=&
\frac{1}{(h-i)!} T_{h-i}(z) \mbox{ for } 0\leq i \leq h \\ \\
S_i(z) &=&
z^{rn} P_{i-h} (z) \mbox{ for } h+1 \leq i \leq \bplush.
 \end{array}
 \right.
\eneq
Then we have for any $q\in\zeroh$:
 
\begin{eqnarray*}
\hh_q(z)
&=& \sum_{p=0}^q \bino (-\log z)^{q-p} \Big( \TT_p(z) +\sum_{i=p+1}^{p+b} \qqp_i(z)g_i(z)\Big) \\
&&\quad \mbox{using Eqns.~\eqref{eqdefff} and~\eqref{eqdefhh}, since $\qqp_i(z)=0$ if $i\leq p$ or $i\geq b+p+1$} \\
&=& \sum_{p=0}^q \bino (-\log z)^{q-p} \TT_p(z) + \sum_{p=0}^q \bino (-\log z)^{q-p} \sum_{i=1}^b z^{rn} P_i(z)
 (-1)^p (i)_p g_{i+p}(z) \\
 &&\quad \mbox{using Eq.~\eqref{eqdefqp}} \\
&=& \sum_{i=h-q}^h \frac{1}{(h-i)!} \TT_{h-i}(z) \frac{q!}{(i+q-h)!} (-\log z)^{i+q-h} \\
 &&\quad \quad + \sum_{i=h+1}^{\bplush} z^{rn} P_{i-h}(z) \sum_{p=0}^q \bino (-\log z)^{q-p} (-1)^p (i-h)_p g_{i-h+p}(z) 
 \end{eqnarray*}
so that 
\begin{equation} \label{eqcclhh}
\hh_q(z) = \sum_{i=0}^{\bplush} \SS_i(z) \yy_{i,q}(z) 
\eneq
by definition of $ \SS_i(z)$ and $ \yy_{i,q}(z) $. The point in writing $\hh_q(z)$ in this way is that the polynomials $S_i(z)$ are independent from $p$ (or $q$).

\bigskip

\noindent {\bf Step 3:} A differential system independent from $p$ (or $q$).

The construction in Step 2 has an important feature: the vectors $Y_q = \tra ( \yy_{0,q},\ldots, \yy_{\bplush,q})$ are solutions of the same differential system, independent from $q$. This is what we shall prove now. 

In precise terms, we claim that for any $q\in\zeroh$ we have:
\begin{equation}\label{eqsysdiff}
\left\{
\begin{array}{rcl}
\yy_{i,q}'(z) &=&
- \frac1z \yy_{i-1,q}(z) \mbox{ for } 1 \leq i \leq \bplush \mbox{ such that } i \neq h+1 \\
\yy_{h+1,q}'(z) &=&
\frac{z+1}{z(1-z)} \yy_{h,q} (z) \\
\yy_{0,q}'(z) &=&
0.
 \end{array}
 \right.
\eneq
We shall check this property now by considering successively various ranges for $i$. If $i=0$, we have $\yy_{0,q}(z) =0$ if $q\leq h-1$ and $\yy_{0,h}(z) =h!$. If $1 \leq i \leq h-q-1$ we have $\yy_{i,q}(z) = \yy_{i-1,q}(z) = 0$. If $i=h-q$ then $\yy_{i,q}(z) =q!$ and $ \yy_{i-1,q}(z) = 0$. In the case where $h-q+1\leq i\leq h$, the derivative of $\yy_{i,q}(z) = \frac{q!}{(i+q-h)!} (-\log z)^{i+q-h} $ is equal to $ - \frac1z \frac{q!}{(i+q-h-1)!} (-\log z)^{i+q-h-1} = 
- \frac1z \yy_{i-1,q} (z) $. When $i=h+1$ the derivative of $\yy_{i,q}(z)$ can be computed as follows:
\begin{eqnarray*}
\yy_{h+1,q}'(z) 
&=&
\sum_{p=0}^q \bino (-1)^p p! \Big( - \frac1z (q-p) (-\log z)^{q-p-1} g_{p+1}(z) + (-\log z)^{q-p} g_{p+1}'(z) \Big)\\
&=&
- \frac1z \Big( \sum_{p=0}^{q-1} \frac{q!}{(q-p-1)!} (-1)^p (-\log z)^{q-p-1} g_{p+1}(z) \\
&&
\quad 
+ \sum_{p=1}^q \frac{q!}{(q-p)!} (-1)^p (-\log z)^{q-p} g_{p}(z) \Big)+ (-\log z)^{q} \cdot \frac{z+1}{z(1-z)} g_0(z) \\
&&
\quad \quad \mbox{since $g_{p+1}'(z) = - \frac1z g_{p}(z) $ for $p\geq 1$, and $g_{1}'(z) = \frac{z+1}{z(1-z)} g_0(z) $}\\
&=&  \frac{z+1}{z(1-z)} \yy_{h,q}(z) 
 \end{eqnarray*}
since the two sums inside the bracket are opposite of each other, and $g_0$ is the constant function equal to $x_0=1$. At last, for $h+2\leq i \leq \bplush$ we have a similar computation:
\begin{eqnarray*}
\yy_{i,q}'(z) 
&=&
- \frac1z \Big( \sum_{p=0}^{q-1} \frac{q!}{(q-p-1)!} (-1)^p \frac{(i-h)_p}{p!} (-\log z)^{q-p-1} g_{i-h+p}(z) \\
 &&\quad \quad + \sum_{p=0}^q \frac{q!}{(q-p)!} (-1)^p \frac{(i-h)_p}{p!} (-\log z)^{q-p} g_{i-h+p-1}(z) \Big) \\
&=&
- \frac1z \sum_{p=0}^{q} \frac{q!}{(q-p)!} (-1)^p (-\log z)^{q-p} g_{i-h+p-1}(z) \Big( - \frac{(i-h)_{p-1}}{(p-1)!} + \frac{(i-h)_p}{p!} \Big)
 \end{eqnarray*}
where $ \frac{(i-h)_{p-1}}{(p-1)!} $ should be understood as 0 for $p=0$. Now $ - \frac{(i-h)_{p-1}}{(p-1)!} + \frac{(i-h)_p}{p!} = \frac{(i-h-1)_{p}}{p!}$ for any $p\geq 0$, so that $\yy_{i,q}'(z) = - \frac1z \yy_{i-1,q} (z)$. This concludes the proof of the claim.

\bigskip

\noindent {\bf Step 4:} Linear independence of the functions $\hh_0$, \ldots, $\hh_h$.

Recall that $\hh_q$ has been defined in Step 1 by Eq.~\eqref{eqdefhh}, for $q\in\zeroh$. Let us prove that these functions are linearly independent over $\C$. Let $\lam_0$, \ldots, $\lam_h\in\C$ be such that $\sum_{q=0}^h \lam_q \hh_q(z)=0$. Then Eq.~\eqref{eqcclhh} yields
\begin{equation}\label{equab}
 \sum_{i=0}^{\bplush} S_i(z) \sum_{q=0}^h \lam_q y_{i,q}(z)=0.
\eneq
Now let $y_i(z)=\sum_{q=0}^h \lam_q y_{i,q}(z)$ for $0\leq i \leq \bplush$. Then Eqns.~\eqref{eqsysdiff} yield $y_0'(z)=0$, $y_{h+1}'(z)= \frac{z+1}{z(1-z)}y_h(z)$, and $y_i'(z)=-\frac1z y_{i-1}(z)$ for any $i \in\{1,\ldots,\bplush\}\setminus\{h+1\}$. 

Assume that $\lam_0$, \ldots, $\lam_h$ are not all zero. Let $q_0$ be the maximal index $q\in\zeroh$ such that $\lam_q\neq0$. Then Eqns.~\eqref{eqdefyq} yield $y_{h-q_0}(z)= \sum_{q=0}^{q_0} \lam_q y_{h-q_0,q}(z) = \lam_{q_0} q_0!\neq 0$ and $y_i(z)=0$ for $0\leq i \leq h-q_0-1$. We write $i_0 = h-q_0$, so that $y_{i_0}(z)= \lam_{q_0} q_0!\neq 0$ and $y_i(z)=0$ for $i<i_0$. 

We shall prove by decreasing induction on $\alp\in\{i_0,\ldots,\bplush\}$ that there exist polynomials $U_{\alp,i_0}$, \ldots, $U_{\alp,\alp}$ such that 
\begin{equation}\label{eqrecdes}
U_{\alp,\alp} \mbox{ is not the zero polynomial and } \sum_{i=i_0}^{\alp} U_{\alp,i}(z) y_{i}(z)=0 \mbox{ for any } z\in D,
\eneq
where $D$ is the open disk we have chosen around $z_0$. This is true for $\alp=\bplush$ by definition of $i_0$, upon letting $U_{\bplush,i}(z)=S_i(z)$: recall that $S_{\bplush}(z)=z^{rn}P_b(z)$ is not the zero polynomial (by definition of $b$ at the beginning of \S \ref{subseczero}), and that~\eqref{equab} holds. Assume that~\eqref{eqrecdes} holds for some $\alp\in \{i_0+1,\ldots,\bplush\}$ and denote by $d$ the degree of $U_{\alp,\alp} $. Then the $(d+1)$-th derivative of the zero function can be written as 
$$z^{d+1} (1-z)^{d+1} \Big( \sum_{i=i_0}^{\alp} U_{\alp,i}(z) y_{i}(z) \Big)^{(d+1)} = \sum_{i=i_0}^{\alp-1} U_{\alp-1,i}(z) y_{i}(z)$$
for some polynomials $U_{\alp-1,i}$, using the expression of $y_i'(z)$ in terms of $y_{i-1}(z)$ deduced above from Eqns.~\eqref{eqsysdiff}; notice that $y_\alp(z)$ does not appear anymore since $U_{\alp,\alp} ^{(d+1)}=0$. 
To prove that $U_{\alp-1,\alp-1}\neq 0$, we first assume that $\alp\neq h+1$. By induction on $t\geq 0$ we have
$$ \Big( U_{\alp,\alp}(z) y_\alp(z)\Big)^{(t)} = U_{\alp,\alp}^{(t)}(z) y_\alp(z) + \sum_{j=0}^{t-1} \Big( \frac{-1}{z} U_{\alp,\alp}^{(j)}(z) \Big) ^{(t-1-j)} y_{\alp-1}(z) + V_t(z)$$ 
where $z^{t} (1-z)^{t} V_t(z)$ is a $\Qbar[z]$-linear combination of $y_{i_0}(z)$, \ldots, $y_{\alp-2}(z) $. Therefore we have 
$$U_{\alp-1,\alp-1}(z) = z^{d+1} (1-z)^{d+1} \Big(U_{\alp,\alp-1}^{(d+1)}(z) + \sum_{j=0}^{d} \Big( \frac{-1}{z} U_{\alp,\alp}^{(j)}(z) \Big) ^{(d-j)} \Big).
$$ 
This is not the zero polynomial because in the expansion of $z^{-(d+1)} (1-z)^{-(d+1)} U_{\alp-1,\alp-1}(z)$ as a linear combination of $z^n$, $n\in\Z$, the coefficient of $z^{-1}$ (namely, the residue) is $- U_{\alp,\alp}^{(d)}\neq 0$. In the case where $\alp = h+1$ we have $y_\alp'(z)=\frac{z+1}{z(1-z)}y_{\alp-1}(z)$ so that the same formulas hold with $ \frac{z+1}{z(1-z)}$ instead of $\frac{-1}{z}$; we conclude in the same way, by writing $z^{-(d+1)} (1-z)^{-(d+1)} U_{\alp-1,\alp-1}(z) = \sum_{n=n_0}^{+\infty} a_n z^n$ for some $n_0\leq -1$, with $a_{-1}\neq 0$. In both cases this concludes the inductive proof of~\eqref{eqrecdes} for all $\alp\in \{i_0,\ldots,\bplush\}$. 

Now for $\alp=i_0$ we obtain $U_{i_0,i_0}(z) y_{i_0}(z)=0$ for any $z\in D$, where $ U_{i_0,i_0} $ is not the zero polynomial and $ y_{i_0}(z)=\lam_{q_0} q_0!\neq 0$. This contradiction concludes the proof of the claim.

\bigskip

\noindent {\bf Step 5:} Defining linearly independent functions $\hhti_1$, \ldots, $\hhti_b$.

Consider, for $\beta\in \unb$, the functions $\yti_{i,\beta}$ defined by
\begin{equation}\label{eqdefyyti}
\left\{\begin{array}{l}
\yti_{i,\beta}(z) = 0 \mbox{ for } 0\leq i \leq h+\beta-1\\
\yti_{i,\beta}(z) = \frac{(-\log z)^{ i-h-\beta }}{ ( i-h-\beta )!} \mbox{ for } h+\beta \leq i \leq \bplush
\end{array}\right.
\eneq
They satisfy the differential system~\eqref{eqsysdiff}; we define
\begin{equation}\label{eqdefhhti}
\hhti_\beta(z) = \sum_{i=0}^{\bplush} S_i(z) \yti_{i,\beta}(z) = \sum_{i=h+\beta}^{\bplush} z^{rn} P_{i-h}(z) \frac{(-\log z)^{ i-h-\beta }}{ ( i-h-\beta )!}= \sum_{i= \beta}^{b} z^{rn} P_{i }(z) \frac{(-\log z)^{ i -\beta }}{ ( i -\beta )!}.
\eneq
Let us prove that the functions $\hhti_1$, \ldots, $\hhti_b$ are linearly independent over $\C$. Let $\lam_1$, \ldots, $\lam_b$ be complex numbers, not all zero, such that $\sum_{\beta=1}^b \lam_\beta \hhti_\beta(z) =0$. Denote by $\beta_0$ the least index $\beta$ such that $ \lam_\beta\neq 0$. Then we have the following $\C[z]$-linear relation between powers of $\log z$:
$$\sum_{\beta=\beta_0}^b \sum_{i=\beta}^b \lam_\beta z^{rn} P_{i }(z) \frac{(-\log z)^{ i -\beta }}{ ( i -\beta )!} =0.$$
Since $\log z$ is transcendental over $\C[z]$, the coefficient of $( \log z)^{ b -\beta_0 }$ has to be zero: $\lam_{\beta_0}P_b(z)=0$. Since $ \lam_{\beta_0}\neq 0$ and $P_b$ is not the zero polynomial (by definition of $b$, see the beginning of \S \ref{subseczero}), this is a contradiction. This concludes the proof that $\hhti_1$, \ldots, $\hhti_b$ are linearly independent over $\C$.

\bigskip

\noindent {\bf Step 6:} Application of Shidlovsky's lemma.

Let us apply the general version of Shidlovsky's lemma stated as 
Theorem~\ref{thzerofct} in \S \ref{subsecshidenonce}.
 We let $\qshid = \bplush+1$ and consider the matrix $A\in M_\qshid(\Q(z))$ that corresponds to the differential system~\eqref{eqsysdiff}. The polynomials $S_0,\ldots,S_{\bplush}$ are defined by Eq.~\eqref{eqdefSS}; we have $\deg S_i \leq m$ with $m=2rn$ (recall that $r\geq 1$, $\deg T_p\leq 2rn$ and $\deg P_i\leq n$). 
 We let $\Sigma = \{0,1,\infty,z_0\}$; recall that $z_0\not\in\{0,1\}$. Let us start with the vanishing conditions at $z_0$.
 
 Eq.~\eqref{eqcclhh} reads $R(Y_q)(z)=\hh_q(z)$ for any $q\in\zeroh$, where $Y_q = \tra (y_{0,q}(z),\ldots,y_{\bplush,q}(z))$ is a solution of $Y'=AY$. The functions $y_{i,q}(z)$ are analytic at $ z_0$ (since $z_0\not\in\{0,1\}$), and the remainders $R(Y_q)(z)=\hh_q(z)$, for $q\in J_{z_0}=\zeroh$, are linearly independent over $\C$ (as proved in Step~4). Moreover we have proved in Step~1 that $\ff_p(z) = O((z-z_0)^{\capa n})$ as $z\to z_0$, so that $R(Y_q)(z)= O((z-z_0)^{\capa n})$ for any $q$ using Eq.~\eqref{eqdefhh}. Therefore we have 
\begin{equation}\label{eqshidzz}
 \sum_{j\in J_{z_0} } \ord_{z_0}(R(Y_j)) \geq (h+1) \capa n.
 \eneq
 
Let us consider now the points 0 and $\infty$. We let $J_0=J_\infty=\unb$, and for $\beta$ in this set we let $\gdyti_\beta = \tra (\yti_{0,\beta}(z),\ldots,\yti_{\bplush,\beta}(z))$ where the functions $\yti_{i,\beta}(z)$ have been defined in Step~5. Then $R(\gdyti_\beta)(z)= \hhti_\beta(z)$ is given by Eq.~\eqref{eqdefhhti}; as proved in Step 5, the functions $R(\gdyti_1)$, \ldots, $R(\gdyti_b)$ are 
$\C$-linearly independent. Recall from Eq.~\eqref{eqdefSS} that $S_i(z)=O(z^{rn})$ as $z\to 0$, and $\deg S_i\leq (r+1)n$, for any $i\in\{h+1,\ldots,\bplush\}$. Therefore Eqns.~\eqref{eqdefyyti} and~\eqref{eqdefhhti} yield $\hhti_\beta(z) =O(z^{rn}(\log z)^{b-1})$ as $z\to 0$, and $\hhti_\beta(z) =O((1/z)^{-(r+1)n}(\log (1/z))^{b-1})$ as $z\to \infty$, so that 
\begin{equation}\label{eqshidzi}
\sum_{\sigma\in\{0,\infty\}} \sum_{\beta\in J_{\sigma} } \ord_{\sigma}(R(\gdyti_\beta)) \geq brn - b(r+1)n = -bn;
 \eneq
 recall that logarithmic factors have no influence on the order of vanishing, e.g. $\ord_0(z^e(\log z)^i)= \mbox{Re}(e)$ for $e\in\C$ and $i\in\N$. 
 
 At last, we let $J_1=\{1\}$ and notice that $R(\gdyti_1)(z)= \hhti_1(z)$ defined by Eq.~\eqref{eqdefhhti} is equal to $z^{rn} R_n(z)$, where $R_n(z)$ is defined in Eq.~\eqref{eqdefrn} (recall that $P_{b+1}(z)=\ldots=P_a(z)=0$). The proof of Theorem~\ref{thsiegel} (namely $(iii)$ in \S \ref{subsecdio2}) shows that $R_n(z) = O((z-1)^{\om n-1})$ as $z\to 1$; therefore we have
 \begin{equation}\label{eqshidun}
 \ord_{1}(R(Y_1)) \geq \om n-1 
 \eneq
where $R(Y_1)$ is not the zero function (see Step 5). 

Combining Eqns.~\eqref{eqshidzz},~\eqref{eqshidzi} and~\eqref{eqshidun}, Theorem~\ref{thzerofct} yields
$$ \Big( (h+1)\capa -b+\om\Big) n - 1 \leq (2rn+1) (\mu-b) + \cstun$$
where $\cstun$ depends only on $a$, $h$, $z_0$ (but can be made independent of $b$ and  $n$ since $b\leq a$), and $\mu$ is the minimal order of a non-zero differential operator $L $ such that $L(R(Y))= 0 $ for any solution $Y$ of the differential system $Y'=AY$. 
Now for any such $Y$, the row matrix $\tra (R(Y)\, R(Y)' \ldots R(Y)^{(N)})$ can be written as $\tra Y M$ where $M\in M_{N,N+1}(\Qbar(z))$ is independent of $Y$: the first column of $M$ is given by the $S_i$, and the following ones by rational functions $S_{k,i}$ (see~\cite[\S 3.2, Step~1]{SFcaract}). There is a non-trivial $\Qbar(z)$-linear relation between the columns of $M$; it provides a differential operator $L$ of order at most $N$ such that $L(R(Y))= 0 $ for any solution $Y$, so that 
 $\mu\leq \bplush+1$. Since $n$ is assumed to be sufficiently large (in terms of $b$, $h$, $\om$, $r$, $z_0$ and $\capa$, and also therefore in terms of $\cstun$), we obtain $
(h+1)(\capa-2r)+\om \leq b $. Since $b \leq a $, $\om>0$ and $(h+1)(\capa-2r)+\om > a$, this is a contradiction.

\subsection{End of the proof} \label{subsecfin}

Let $a$ be a positive integer. In Theorem~\ref{thintroun} the numerical constant $0.21$ can be replaced (as the proof will show) by a slightly larger real number. Therefore in the proof we may assume that $a$ is a multiple of 25. Then we choose $r =3.9$, $\capa = 10.58$, $\om = 12$, $\gdom = \lfloor r \sqrt{a \log a}\rfloor$, and $h =0.36\ a\in\N$, so that $(h+1)(\capa-2r)+\om > a$ and $\gdom> \om$. Here and below all numerical constants are rounded with precision 0.01.

We consider $z_0 = -1 $ and choose $\qq=1$, so that $\qq z_0 \in \Z$. We denote by $\calN_a$ the set of all sufficiently large multiples of 50: 
for any $n\in\calN_a$ we have $rn , \capa n , \om n,\gdom n\in \N$. For any $n\in\calN_a$ we 
consider the integers $\cij $ provided by Theorem~\ref{thsiegel}, and define $b$ as in the beginning of \S \ref{subseczero}, namely
$$b = \max\{ i\in\{1,\ldots,a\}, \, \exists j \in\zeron, \, \cij\neq 0\}.$$
This integer $b$ depends on $n$, but it can take only $a$ values. Therefore  there exists an infinite subset $\calN_a'\subset\calN_a$ such that all $n\in\calN_a' $ correspond to the same $b$. From now on, we consider only integers $n\in\calN_a'$.

Let  $k\in\{2rn+2,\ldots,\capa n\}$ and $p \in\{ 0 ,\ldots,h\}$.  Lemma~\ref{lemcoeffs} yields  $\likn\in \Z$ for any $i$,   and 
$$|  \likn  |  \leq \beta^{n(1+o(1))} \mbox{ with } \beta =
 \chi \Big( 8 e^3(2a+1)\Big)^\capa \cdot 2^{\capa+r+1}
$$
where $\chi$ is defined by Eq.~\eqref{eqdefchi} in Theorem~\ref{thsiegel}, namely
$$
\chi = \exp\Big(\frac{ \om \log 2 + 3\om^2 + \om^2 \log (a+1) +\frac12 \gdom^2 \log r }{a-\om }\Big).
$$
Now we have (using Eq.~\eqref{eqFL} and the definition of $b$, see the beginning of \S \ref{subseczero})
$$ \lzkn + \sum_{i=1}^{\bplush} \likn \Big( 1 - (-1)^i\Big) \Li_i(-1) = (-2)^{k-1} \frac{\del_{k}}{(k-1)!}
S_{n,p}^{(k-1)}(-1).$$
Since $k \leq\capa n$, we may apply Lemma~\ref{lemasy} and deduce that 
$$\Big| \lzkn + \sum_{i=1}^{\bplush} \likn \Big( 1 - (-1)^i\Big) \Li_i(-1) \Big| \leq \alpha^{n(1+o(1))} \mbox{ with } 
\alpha =
2^\capa \alpha_0 = \chi r^{-\gdom}( 2e^4(2a+1))^{\capa} . $$
Using Proposition~\ref{proplemzero}, the refined version of  Siegel's linear independence criterion (stated and proved in  \S \ref{subsecsiegel})
applies to these linear forms  for $n\in \calN_a' $, with coefficients $\likn$, $\theta_0=1$, $Q_n = \beta^n$ and $\tau = -\frac{\log \alpha}{\log\beta}$ (so that $Q_n^{-\tau} = \alpha^n$). We obtain 
\begin{equation} \label{eqaveclogd}
\dim_\Q \Span_\Q ( \{1, \log 2\}\cup\{\zeta(i), \, 3\leq i \leq a+h, \, i \mbox{ odd}\}) \geq 1 - \frac{\log \alpha}{\log \beta } .
\eneq
Now recall that $a>0$ is a multiple of $25$, 
$r =3.9$, $\capa = 10.58$, $\om = 12$, $\gdom =\lfloor r \sqrt{a \log a}\rfloor$, and $h =0.36\ a$.
 As $a\to\infty$ the formulas above yield
 $$\log\chi\sim \frac{ \gdom^2 \log r}{2a} \sim \frac{r^2\log r}2 \log a ,$$
$$\log\beta\sim \log \chi + \capa \log a \sim \Big( \frac{r^2\log r}2 +\capa\Big) \log a , $$ 
$$\log\alpha\sim - \gdom \log r \sim - r\log r \cdot \sqrt{a \log a} $$
so that 
\begin{eqnarray*}
- \frac{\log \alpha}{\log \beta} 
&\sim& \frac{ 2 r \log r}{r^2 \log r + 2 \capa} \sqrt{\frac{a}{\log a}} \\
&\sim& \frac{ 2 r \log r}{r^2 \log r + 2 \capa}\cdot \frac1{\sqrt{1+h/a}} \cdot \sqrt{\frac{a+h}{\log (a+h)}}.
\end{eqnarray*}
Now recall that $r=3.9$, $\capa=10.58$ and $h=0.36 a$, so that 
$$ \frac{ 2 r \log r}{r^2 \log r + 2 \capa}\cdot \frac1{\sqrt{1+h/a}} = 0.2174\ldots > 0.21.$$
If $a$ is large enough we obtain
$$- \frac{\log \alpha}{\log \beta} \geq 0.21 \cdot \sqrt{\frac{a+h}{\log (a+h)}}.$$
We take $s = a+h$ and apply Eq.~\eqref{eqaveclogd}. The additional 1 in the right hand side accounts for the number $\log 2$ in the left hand side, that we want to get rid of.
This concludes the proof of Theorem~\ref{thintroun}.

\bigskip

\begin{Remark} It follows from the computations above that, as $s = a+h$ tends to $\infty$, 
$$
\log\alpha\sim - 4.55 \sqrt{s \log s} \quad \mbox{ and } \quad \log\beta\sim 20.93 \log s.
$$

\end{Remark}

\begin{Remark}
The proof allows one to compute effectively an integer $s_0$ such that the conclusion of Theorem~\ref{thintroun} holds for any $s\geq s_0$.
\end{Remark}

\subsection{The case of polylogarithms: proof of Theorem~\ref{thpolylogs}} \label{subsecpolylogs}

To prove Theorem~\ref{thpolylogs}, we follow the proof of Theorem~\ref{thintroun} except that we consider $\Snpinf(z)$
(defined in Eq.~\eqref{eqdefsninf})
 instead of $S_{n,p}(z)$. Therefore Eq.~\eqref{eq33nv} becomes
\begin{equation} \label{eq33polylogsnv}
 \Snpinfkmu= \qqp_{k,0}(z) + \sum_{i=1}^{\aplush} \qqp_{k,i}(z) \Li_i(1/z) \mbox{ for any } k\geq (r+1)n+1.
\eneq
 The point here is that (with the notation of the proof of Lemma~\ref{lemdev} in \S \ref{subsecdev})
we have $\deg \Vpinf \leq (r+1)n-1$ and $\deg \Vpzero \leq 2rn$. In the proof of Theorem~\ref{thintroun} we had to restrict to integers $k\geq 2rn+2$ so that $ ( \Vpinf - \Vpzero) ^{(k-1)}=0$, whereas to prove Theorem~\ref{thpolylogs} assuming $ k\geq (r+1)n+1$ is enough to ensure that $ \Vpinfkmu =0$. 

Let $z_0\in\Qbar$ be such that $|z_0| \geq 1$ and $z_0\neq 1$; let $\qq\in\N\etoile$ be a denominator of $z_0$, i.e. such that $\qq z_0 \in \OQzz$ where $\OQzz$ is the ring of integers of $\Q(z_0)$.
For any $k\geq (r+1)n+1$ we let 
$$
\likn(z_0) = q^{(r+1)n+k-1} z_0^{k-1}(1-z_0)^{k-1} \frac{\del_k}{(k-1)!} \qqp_{k,i} (z_0) \mbox{ for } 0\leq i \leq \aplush 
$$
where $\del_k=\del_k (a+h,(r+1)n, 1,0)$ is given by Proposition~\ref{propgene} in \S \ref{subsectech1} with $a $ replaced by $\aplush$ and $n$ by $(r+1)n$; 
in the setting of \S \ref{subsectech1} we take $\bte = 1$ and $ \ate = 0 $ in the recurrence relation~\eqref{eqrecpkijtechnique}, to fit the differential system satisfied by the functions 1 and $ \Li_i(1/z)$.
Then following the proof of Lemma~\ref{lemcoeffs} (with only one difference: for $i=0$, due to the value of $(\bte,\ate)$) yields $ \likn(z_0)\in\OQzz$ and 
$$ \house{\likn(z_0)} \leq \beta_1^{n(1+o(1))} \mbox{ with } \beta_1 = \chi \Big( 8e^3 (2a+1) \Big)^\capa \cdot \Big( \qq \max(1, \house{z_0} , \house{ 1-z_0} ) \Big)^{\capa+r+1} $$
provided $k \leq \capa n$ and $\capa\geq r+1$. Moreover Eq.~\eqref{eq33polylogsnv} yields 
$$ q^{(r+1)n+k-1} z_0^{k-1}(1-z_0)^{k-1} \frac{\del_k}{(k-1)!} \Snpinfkmu(z_0) = \lzkn(z_0) + \sum_{i=1}^{\aplush} \likn(z_0) \Li_i(1/z_0) 
$$
for any $ k\geq (r+1)n+1$. Following the proof of Lemma~\ref{lemasy} we deduce that 
$$\Big| q^{(r+1)n+k-1} z_0^{k-1}(1-z_0)^{k-1} \frac{\del_k}{(k-1)!} \Snpinfkmu(z_0) \Big| \leq \alpha_1^{n(1+o(1))} $$
 with 
$$\alpha_1 = \chi r^{-\gdom} q^{r+1} ( e^4(2a+1)q | z_0(1-z_0)| )^{\capa} .$$
Then we adapt Proposition~\ref{proplemzero}, assuming that $(h+1)(\capa- r-1)+\om > a$ and considering integers $k$ such that $(r+1)n+1\leq k \leq \capa n$. This enables us to apply Proposition \ref{propsiegel}
 and deduce that 
$$\dim_{\Q(z_0)} \Span_{\Q(z_0)} (\{1\}\cup \{ \Li_i(1/z_0), \, 1\leq i \leq \aplush\}) \geq \frac1{[ \Q(z_0):\Q]} \Big( 1 - \frac{\log \alpha_1}{\log \beta_1 } \Big).$$
Our choice of parameters is the same as in \S \ref{subsecfin}, except for numerical constants. The only difference is that the assumptions $ \capa> 2 r$ and $(h+1)(\capa-2 r)+\om > a$ in \S \ref{subsecfin} are weakened here to $ \capa> r+1$ and $(h+1)(\capa- r-1)+\om > a$. We choose 
$r =5.3$, $\capa = 8.8343$, $\om = 10$, $\gdom = \lfloor 3.3 \sqrt{a \log a}\rfloor$, and $h =0.3946\ a\in\N$ (assuming that $10^4$ divides $a$), so that $(h+1)(\capa-r-1)+\om > a$. As in \S \ref{subsecfin} we have, as $a\to\infty$:
$$\log\chi\sim 9.0807 \log a, \quad 
 \log\beta_1\sim 17.915 \log a , \quad 
\log\alpha_1\sim - 5.5034 \sqrt{a \log a} $$
so that 
$$ - \frac{\log \alpha_1}{\log \beta_1} \geq 0.26 \sqrt{\frac{a+h}{\log (a+h)}}
$$
provided $a$ is large enough. This concludes the proof of Theorem~\ref{thpolylogs}.

\begin{Remark} If $z\not\in\R$ then we have $[\K_\infty:\R]=2$ in the notation of Proposition~\ref{propsiegel}, so that the constant $0.26$ may be replaced with $0.52$ in Theorem~\ref{thpolylogs}.
\end{Remark}

\providecommand{\bysame}{\leavevmode ---\ }
\providecommand{\og}{``}
\providecommand{\fg}{''}
\providecommand{\smfandname}{\&}
\providecommand{\smfedsname}{\'eds.}
\providecommand{\smfedname}{\'ed.}
\providecommand{\smfmastersthesisname}{M\'emoire}
\providecommand{\smfphdthesisname}{Th\`ese}

\end{document}